\documentclass[a4paper, reqno, 11pt]{amsart}

\usepackage{geometry} 
\usepackage{latexsym}
\usepackage{fullpage}
\usepackage{amssymb}
\usepackage{amscd}
\usepackage{float}
\usepackage{graphicx}
\usepackage{amsfonts}
\usepackage{pb-diagram}
\usepackage{amsmath,amscd}
\usepackage{caption}
\usepackage{subcaption}
\usepackage{xcolor}

\newtheorem{theorem}{Theorem}

\newtheorem{definition}{Definition}
\newtheorem{remark}{Remark}
\newtheorem{example}{Example}

\UseRawInputEncoding

\begin{document}
	
\title[A survey on skein modules via braids]
  {A survey on skein modules via braids}

\author{Ioannis Diamantis}
\address{Department of Data Analytics and Digitalisation,
Maastricht University, School of Business and Economics,
P.O.Box 616, 6200 MD, Maastricht,
The Netherlands.
}
\email{i.diamantis@maastrichtuniversity.nl}



\date{}

\begin{abstract}
In this paper we present recent results on the computation of skein modules of 3-manifolds using braids and appropriate knot algebras. Skein modules generalize knot polynomials in $S^3$ to knot polynomials in arbitrary 3-manifolds and they have become extremely essential algebraic tools in the study of 3-manifolds. In this paper we present the braid approach to the HOMFLYPT and the Kauffman bracket skein modules of the Solid Torus ST and the lens spaces $L(p,1)$ and $S^1\times S^2$. 
\end{abstract}

\maketitle

\section{Introduction}\label{intro}

Skein modules were introduced independently by Przytycki \cite{P} and Turaev \cite{Tu} as generalizations of knot polynomials in $S^3$ to knot polynomials in arbitrary 3-manifolds. They are quotients of free modules over isotopy classes of links in 3-manifolds by properly chosen local (skein) relations. A skein module of a 3-manifold yields all possible isotopy invariants of knots which satisfy a particular skein relation. Skein modules based on the Kauffman bracket skein relation 
\[
L_+-AL_{0}-A^{-1}L_{\infty}
\] 
\noindent where $L_{\infty}$ and $L_{0}$ are represented schematically by the illustrations in Figure~\ref{skein}, are called {\it Kauffman bracket skein modules} (KBSM). 

\smallbreak

The precise definition of KBSM is as follows:

\begin{definition}\rm
Let $M$ be an oriented $3$-manifold and $\mathcal{L}_{{\rm fr}}$ be the set of isotopy classes of unoriented framed links in $M$. Let $R=\mathbb{Z}[A^{\pm1}]$ be the Laurent polynomials in $A$ and let $R\mathcal{L}_{{\rm fr}}$ be the free $R$-module generated by $\mathcal{L}_{{\rm fr}}$. Let $\mathcal{S}$ be the ideal generated by the skein expressions $L_+-AL_{0}-A^{-1}L_{\infty}$ and $L \bigsqcup {\rm O} - (-A^2-A^{-2})L$. Note that blackboard framing is assumed and that $L \bigsqcup {\rm O}$ stands for the union of a link $L$ and the trivially framed unknot in a ball disjoint from $L$. 

\begin{figure}[H]
\begin{center}
\includegraphics[width=2.1in]{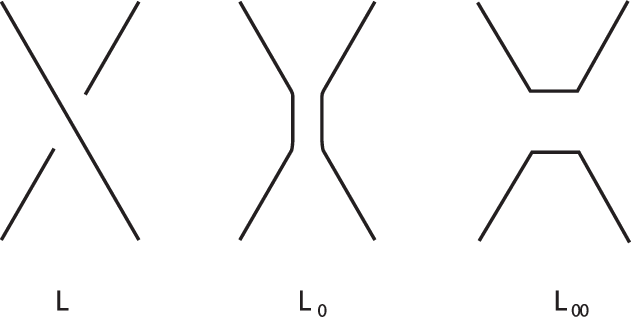}
\end{center}
\caption{The links $L$, $L_{0}$ and $L_{\infty}$ locally.}
\label{skein}
\end{figure}

\noindent Then the {\it Kauffman bracket skein module} of $M$, KBSM$(M)$, is defined to be:

\begin{equation*}
{\rm KBSM} \left(M\right)={\raise0.7ex\hbox{$
R\mathcal{L}_{{\rm fr}} $}\!\mathord{\left/ {\vphantom {R\mathcal{L_{{\rm fr}}} {\mathcal{S} }}} \right. \kern-\nulldelimiterspace}\!\lower0.7ex\hbox{$ S  $}}.
\end{equation*}

\end{definition}

For example, the Kauffman bracket skein module of $S^3$, KBSM($S^3$), is freely generated by the unknot and it is equivalent to the Kauffman bracket for framed links in $S^3$ up to regular isotopy. It's ambient isotopy counterpart is the well-known Jones polynomial. Recall that the algebraic counterpart construction of the Kauffman bracket is the Temperley-Lieb algebra together with a unique Markov trace constructed by V.F.R. Jones (see \S~\ref{Jpoly} or \cite{Jo} and references therein). It is worth pointing out the importance of the Temperley-Lieb algebras. Through the pioneering work of V.F.R. Jones (\cite{Jo, Jo1}), these algebras related knot theory to statistical mechanics, topological quantum field theories and the construction of quantum invariants for 3-manifolds (works of Witten, Reshetikhin-Turaev, Lickorish, etc).

\bigbreak

\begin{figure}[!ht]
\begin{center}
\includegraphics[width=1.7in]{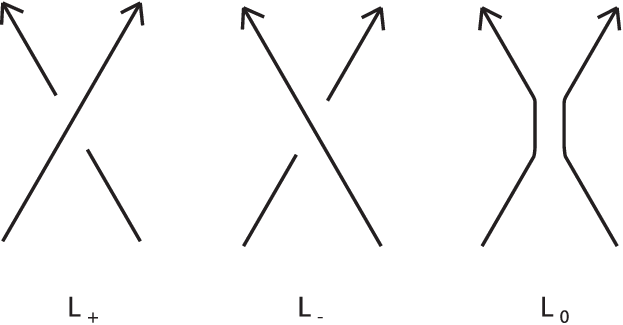}
\end{center}
\caption{The links $L_{+}, L_{-}, L_{0}$ locally.}
\label{skein1}
\end{figure}

Skein modules based on the HOMFLYPT skein relation 
\[
u^{-1}L_{+}-uL_{-}-zL_{0}
\] 
\noindent where $L_{+}$, $L_{-}$ and $L_{0}$ are represented schematically by the illustrations in Figure~\ref{skein1}, are called {\it HOMFLYPT skein modules} (HOM). 

\smallbreak

We have the following definition:

\begin{definition}\rm
Let $M$ be an oriented $3$-manifold, $R=\mathbb{Z}[u^{\pm1},z^{\pm1}]$, $\mathcal{L}$ the set of all oriented links in $M$ up to ambient isotopy in $M$ and let $S$ be the submodule of $R\mathcal{L}$ generated by the skein expressions $u^{-1}L_{+}-uL_{-}-zL_{0}$, where
$L_{+}$, $L_{-}$ and $L_{0}$ comprise a Conway triple represented schematically by the illustrations in Figure~\ref{skein}.

\noindent For convenience we allow the empty knot, $\emptyset$, and add the relation $u^{-1} \emptyset -u\emptyset =zT_{1}$, where $T_{1}$ denotes the trivial knot. Then the {\it HOMFLYPT skein module} of $M$ is defined to be:

\begin{equation*}
\mathcal{S} \left(M\right)=\mathcal{S} \left(M;{\mathbb Z}\left[u^{\pm 1} ,z^{\pm 1} \right],u^{-1} L_{+} -uL_{-} -zL{}_{0} \right)={\raise0.7ex\hbox{$
R\mathcal{L} $}\!\mathord{\left/ {\vphantom {R\mathcal{L} S }} \right. \kern-\nulldelimiterspace}\!\lower0.7ex\hbox{$ S  $}}.
\end{equation*}

\end{definition}

For example, $\mathcal{S}(S^3)$ is freely generated by the unknot (\cite{FYHLMO, PT}). For a survey on skein modules see \cite{HP2}. It is worth mentioning that the HOMFLYPT skein module of a 3-manifold is very hard to compute. 

\bigbreak

Skein modules of $3$-manifolds have become very important algebraic tools in the study of $3$-manifolds, since their properties renders topological information about the $3$-manifolds. For example, the existence of non-separating spheres and tori in a 3-manifold $M$ is captured by the torsion part of KBSM($M$).

\bigbreak

The paper is organized as follows: In \S~\ref{poly} we recall Jones's construction of the Jones and the HOMFLYPT polynomials via a trace function on appropriate knot algebras. For the sake of completion, in \S~\ref{prelim} we also present the Kauffman bracket polynomial following \cite{K}. In \S~\ref{skst} we recall the setting and results from \cite{La1, La2, LR1, LR2, D0} and we present the universal invariants for knots and links in the Solid Torus ST, of the HOMFLYPT type (\S~\ref{homstpoly}) and of the Kauffman bracket type (\S~\ref{kbsmst}). In \S~\ref{smlp} we generalize these universal invariants for knots and links in the lens spaces $L(p,1)$, taking into consideration the extra isotopy move for knots in $L(p,1)$. This is possible with the use of the new bases of the corresponding skein modules presented in \S~\ref{newbasishomst} for HOM(ST) and in \S~\ref{nbas} for KBSM(ST). Finally, in \S~\ref{kbsms1s2} we present recent results on KBSM($S^1\times S^2$) via braids from \cite{D10}.

\section{Invariants of knots in $S^3$}\label{poly}

\subsection{Preliminaries}\label{prelim}

\begin{figure}[H]
\begin{center}\includegraphics[width=3.4in]{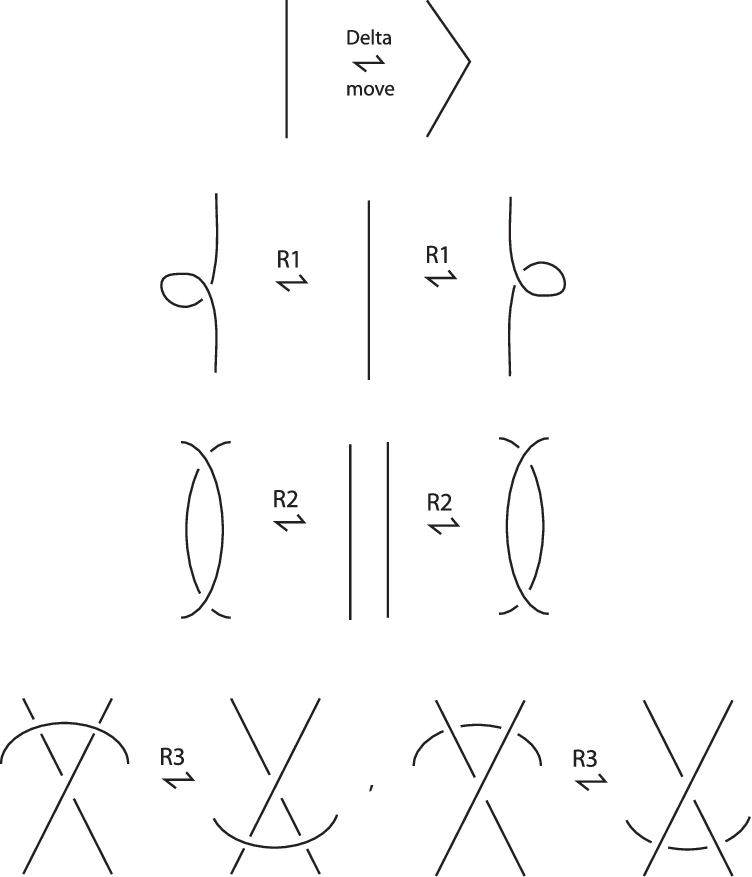}
\end{center}
\caption{Isotopy moves.}
\label{pr1}
\end{figure}

An oriented link on $m$-components is an embedding of $m$-copies of $S^1$ to $S^3$, in which each component is assigned an orientation. The main question of knot theory is {\it ``How can we distinguish one link from another?''} In order to answer this question, we first define an equivalence relation on the set of links. We shall call two links $L_1, L_2$ in $S^3$ equivalent (or isotopic), denoted by $L_1 \sim L_2$, if there is an orientation-preserving piecewise linear homeomorphism $h:S^3\rightarrow S^3$, such that $h(L_1)=L_2$. We may reduce the complexity of the theory by considering link projections on the plane, that we call {\it link diagrams}. More precisely, a link diagram is a projection of the link onto the plane, where we keep track of over and under crossings. Note that we only allow double points and no more than two points are allowed to be superposed. Since now a link diagram depends on the plane that we project the link on, we translate the notion of equivalence of links on the diagrammatic level as follows:

\begin{theorem}[Reidemeister]
Two link diagrams correspond to isotopic links if and only if one can be obtained from the other by a finite sequence of {\it Reidemeister moves} and planar isotopies (or {\it Delta moves}), as illustrated in Figure~\ref{pr1}.
\end{theorem}

\begin{figure}[H]
\begin{center}\includegraphics[width=0.9in]{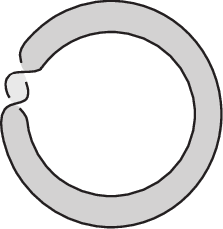}
\end{center}
\caption{A framed knot.}
\label{frk1}
\end{figure}

A framed link on $m$-components is an embedding of a disjoint union of annuli in $S^3$ (see Figure~\ref{frk1}). In the theory of framed links
the Reidemeister move 1 is forbidden.

If two knots are isotopic, using Reidemeister moves we can turn the diagram of one to the diagram of the other. But how can we tell if two knots are not isotopic? We do that by defining functions from the set of links, that we call {\it link invariants}, with the use of which we may determine whether two links are not equivalent. A very important invariant is the Kauffman bracket defined as follows:

\begin{definition}\label{pkaufb}\rm
Let $L$ be a framed link. The {\it Kauffman bracket polynomial} of $L$ is defined by means of the following relations: 
\begin{figure}[H]
\begin{center}
\includegraphics[width=2.3in]{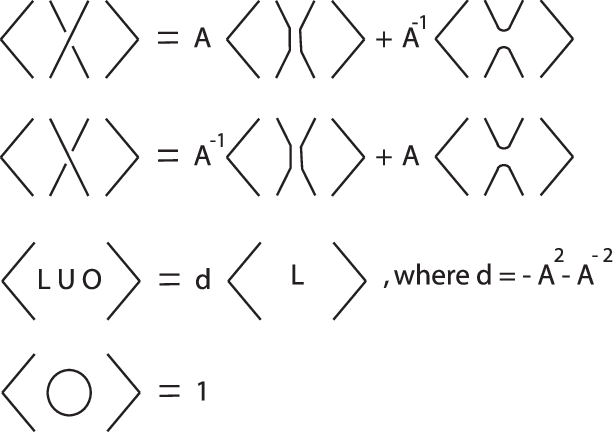}
\end{center}
\label{pkb}
\end{figure}
\end{definition}

\begin{figure}[H]
\begin{center}
\includegraphics[width=1.1in]{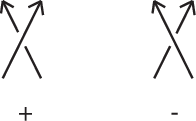}
\end{center}
\caption{The sign of the crossings.}
\label{si}
\end{figure}

Simple computations show that $<L>$, which is a Laurent polynomial in $\mathbb{Z}\left[A^{\pm 1}\right]$, is an invariant of framed links. In order to obtain an invariant for links in $S^3$, we normalize the Kauffman bracket by considering the product of $<L>$ by the factor $\left( -A^3\right) ^{-wr(L)}$, where $wr(L)$ is the writhe of the link $L$, defined as the number of positive crossings minus the number of negative crossings of $L$ (see Figure~\ref{si}).

Then, the polynomial
\[
P_L(A, A^{-1})\ =\ (-A^{-3})^{w(L)}\, <L>,
\]
\noindent is an invariant of links in $S^3$.

\bigbreak

\begin{figure}[H]
\begin{center}\includegraphics[width=1.5in]{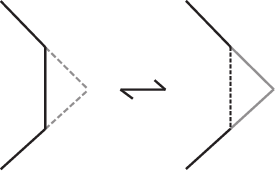}
\end{center}
\caption{Elementary deformations.}
\label{eldef1}
\end{figure}

We now perform isotopy moves on a link diagram and we transform it to a very specific ``knotted object'' that we call a braid. A \textit{braid} on $n$ strands is defined as a set of pairwise non-intersecting descending polygonal lines (strands) joining the points $A_1, A_2, \ldots, A_n$ to the points $B_1, B_2, \ldots, B_n$ in any order, where $A_i=(i,0,0)$ and $B_i=(i,0,1)$ for $i=1, 2, \ldots, n$. Two braids are called \textit{isotopic} if and only if one can be transformed into the other by a finite sequence of {\it elementary deformations} (see Figure~\ref{eldef1}).

\begin{figure}[H]
\begin{center}\includegraphics[width=2.4in]{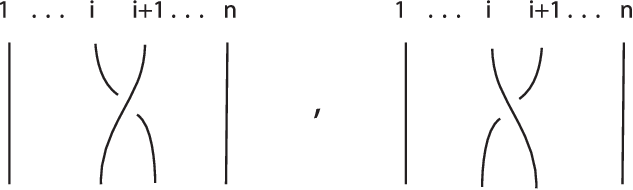}
\end{center}
\caption{The generators $\sigma_i$ and $\sigma_i^{-1}$ of $B_n$.}
\label{br1}
\end{figure}

The set of (equivalent classes of) braids on $n$ strands has a natural group structure. We define the product of two braids $a$ and $b$ by joining the bottom points of the first braid to the top points of the second. The set of braids in $n$ strands under this operation is called the braid group and is denoted by $B_n$. The unit element is the braid consisting of $n$ parallel vertical strands and the inverse of a braid $a$, $a^{-1}$, is the mirror image of $a$ in the plane. $B_n$ has the following presentation:

\[
B_{n} = \left< \begin{array}{ll}  \begin{array}{l} \sigma_{1}, \ldots ,\sigma_{n-1}  \\ \end{array} & \left| \begin{array}{l}
{\sigma_i}\sigma_{i+1}{\sigma_i}=\sigma_{i+1}{\sigma_i}\sigma_{i+1}, \quad{ 1 \leq i \leq n-2}   \\
 {\sigma_i}{\sigma_j}={\sigma_j}{\sigma_i}, \quad{|i-j|>1}  \\
\end{array} \right.  \end{array} \right>,
\]

\noindent where the generators $\sigma_i$ and $\sigma_i^{-1}$ are illustrated in Figure~\ref{br1}.

\begin{figure}[H]
\begin{center}\includegraphics[width=2.1in]{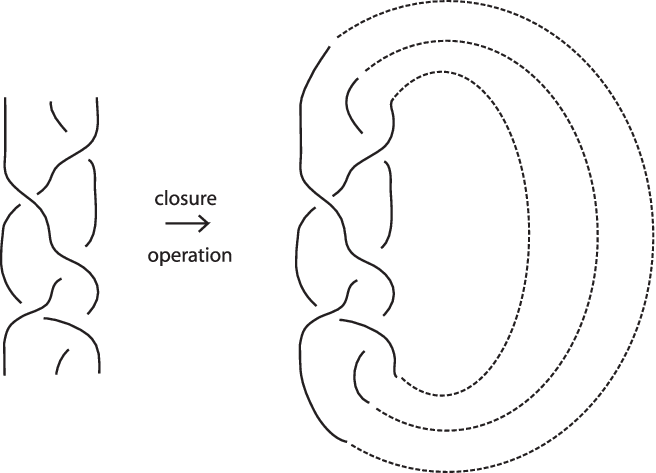}
\end{center}
\caption{The closure of a braid.}
\label{clbr1}
\end{figure}

From a topological point of view, closing a braid, that is, connecting corresponding ends in pairs, gives rise to an oriented link. The \textit{closure of a braid} $a$ is defined as the link $\widehat{a}$ obtained by joining the upper points of its strands to the lower ones (see Figure~\ref{clbr1}). The converse is also true. In particular, we have the following theorem by Alexander:

\begin{theorem}[Alexander]\label{Alex}
Any oriented link is isotopic to the closure of a braid.
\end{theorem}

\begin{figure}[H]
\begin{center}
\includegraphics[width=5.1in]{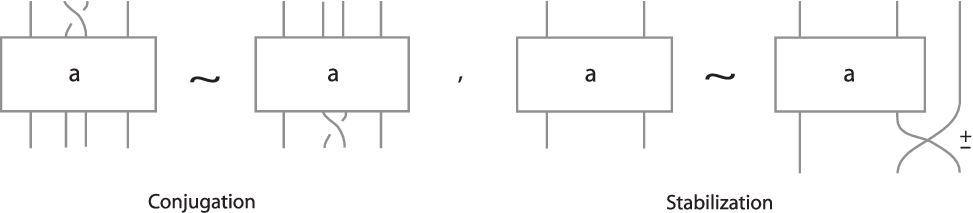}
\end{center}
\caption{Braid equivalence.}
\label{beq}
\end{figure}

Moreover, equivalence classes of braids are in one-to-one correspondence with the isotopy classes of oriented links. More precisely, we have the following theorem:

\begin{theorem}[Markov]\label{Mark}
The closures of two braids are isotopic if and only if one braid can be taken to another by finite sequence of the following moves (for an illustration see Figure~\ref{beq}):
\[
\begin{array}{lllll}
{\rm Conjugation:} & a & \sim & bab^{-1}, & a, b \in B_n, \\
&&&\\
{\rm Stabilization:} & a  & \sim & a\sigma_n^{-1}, & a \in B_n. \\
\end{array}
\]
\end{theorem}

\begin{remark}\rm
Jones noticed that instead of studying knots up to isotopy, we may study braids up to Markov's theorem. That is, 
\[
\frac{{\rm Knots}}{{\rm Reidemeister\ moves}}\ =\ \frac{{\rm Braids}}{\rm Markov's\ theorem}.
\]
This idea lead to the discovery of the Jones polynomial in the 80's.
\end{remark}

\subsection{The HOMFLYPT polynomial}\label{hompoly}

In this subsection we recall the construction of the HOMFLYPT polynomial via a unique Markov trace constructed on the Iwahori-Hecke algebras of type A. The Iwahori-Hecke algebra of type $A$, $H_n$, is a quotient of the braid group algebra $\mathbb{Z}\,  [q^{\pm 1}]B_n$ by factoring out the quadratic relations
\[
\sigma_i^2=(q-1)\sigma_i+q.
\]
\noindent These relations accurately reflect the skein relation of the HOMFLYPT polynomial. The algebra ${\rm H}_n$ has the following presentation:

\[ {\rm H}_{n}(q) = \left< \begin{array}{ll}  \begin{array}{l} g_{1}, g_2, \ldots , g_{n-1}  \\ \end{array} & \left| \begin{array}{l}
g_ig_{i+1}g_i=g_{i+1} g_i g_{i+1}, \quad{ 1 \leq i \leq n-2}   \\
{g_i}{g_j}={g_j}{g_i}, \quad{|i-j|>1}  \\
g_i^2=(q-1)g_i+q, \quad{i=1, 2, \ldots n-1}
\end{array} \right.  \end{array} \right>, \]

\noindent that is

\begin{equation*}
\textrm{H}_{n}(q)= \frac{{\mathbb Z}\left[q^{\pm 1} \right]B_{n}}{ \langle \sigma_i^2 -\left(q-1\right)\sigma_i-q \rangle}.
\end{equation*}

In \cite{Jo} V.F.R. Jones gives the following linear basis for the Iwahori-Hecke algebra of type A, $\textrm{H}_{n}(q)$:

\[
S =\left\{(g_{i_{1} }g_{i_{1}-1}\ldots g_{i_{1}-k_{1}})(g_{i_{2} }g_{i_{2}-1 }\ldots g_{i_{2}-k_{2}})\ldots (g_{i_{p} }g_{i_{p}-1 }\ldots g_{i_{p}-k_{p}})\right\}, \mbox{ for } 1\le i_{1}<\ldots <i_{p} \le n-1{\rm \; }.
\]

\smallbreak

\noindent The basis $S$ yields directly an inductive basis for $\textrm{H}_{n}(q)$, which is used in the construction of the Ocneanu trace, leading to the HOMFLYPT or $2$-variable Jones polynomial. Recall that a linear function $f$ from an algebra to some module is called a
trace if it satisfies $f(xy) = f(yx)$ for all $x, y$ in the algebra. We now have the following:

\begin{theorem}[Ocneanu] \label{ocn}
There exists a unique linear Markov trace function:
\begin{equation*}
{\rm tr}:\bigcup _{n=1}^{\infty }{\rm H}_{n}(q)  \to \mathbb{C},
\end{equation*}

\noindent determined by the rules:

\[
\begin{array}{lllll}
(1) & {\rm tr}(1) & = & 1 &  \\
(2) & {\rm tr}(ab) & = & {\rm tr}(ba) & \quad {\rm for}\ \textit{a,b}\in {\rm H}_n(q) \\
(3) & {\rm tr}(ag_{n}) & = & z{\rm tr}(a) & \quad {\rm for}\ \textit{a}\in {\rm H}_n(q) \\
\end{array}
\]
\end{theorem}

\begin{theorem}\label{hom}
The function $X:\mathcal{L}$ $\rightarrow \mathbb{Z}[q^{\pm1}, z, \lambda]$

\begin{equation*}
X_{L}(q, \lambda)=\left[-\frac{1-\lambda q}{\sqrt{\lambda } \left(1-q\right)} \right]^{n-1} \left(\sqrt{\lambda } \right)^{e}
{\rm tr}\left(\pi \left(\alpha \right)\right),
\end{equation*}

\noindent where $\alpha \in B_{n}$ is a word in the $\sigma _{i}$'s, $e$ is the exponent sum of the $\sigma _{i}$'s in $\alpha $, and
$\pi$ the canonical map of $B_{n}$ in ${\rm H}_{n}(q)$, such that $\sigma _{i} \mapsto g_{i} $, is an invariant of oriented links in ST.
\end{theorem}

\subsection{The Jones polynomial}\label{Jpoly}

We now present the definition of the Temperley-Lieb algebra as a quotient of the Iwahori-Hecke algebra of type $A$, given by Jones \cite{Jo1}.

\begin{definition}\rm
The Temperley-Lieb algebra $TL_n(q)$ is defined as the quotient ${\rm H}_n(q)$ over the ideal $I_n$, where $I_n$ is generated by the element $g_{1, 2}$, where
\[
g_{1, 2}\, :=\, 1\, +\, q\, (g_1\, +\, g_2)\, +\, q^2\, (g_1 g_2\, +\, g_2 g_1)\, +\, q^3\, g_1g_2g_1.
\]
\end{definition}

\begin{remark}\rm
It is worth mentioning that the Temperley-Lieb algebra of type A was originally defined by generators $f_1, \ldots, f_{n-1}$ such that:
\[
\begin{array}{lcll}
f_i^2 & = & f_i &\\
f_i f_j f_i & = & \delta \, f_i & |i-j|=1\\
f_i f_j & = & f_j f_i & |i-j|>1
\end{array}
\]
\noindent where $\delta$ is an indeterminate. If we introduce the generators $g_i:= (q+1)f_i-1$, where $\delta^{-1}= 2 + u + u^{-1}$, then the Temperley-Lieb algebra of type A can be defined as the quotient of $H_n$ over $g_{1, 2}$.
\end{remark}

By annihilating the generator of the defining ideal of the quotient algebra $TL_n$, i.e. by solving $tr(g_{1, 2})=0$, Jones redefined Ocneanu's trace on $TL_n$. In particular, setting $\lambda=q$ in the invariant $X$ of Theorem~\ref{hom}, results in the Jones polynomial $V$. That is, $V(q)\, =\, X(q, q)$.

\begin{remark}\rm
It is important to recall that the Jones polynomial is equivalent to the Kauffman bracket defined in Definition~\ref{pkaufb} up to a change of variable. The Kauffman bracket skein module of a 3-manifold $M$ is equivalent to constructing the Jones polynomial for knots and links in $M$ via the Kauffman bracket.
\end{remark}

\section{Skein modules of the Solid Torus}\label{skst}

In this section we present two universal invariants for knots and links in the Solid Torus ST. One is of the HOMFLYPT type, constructed in \cite{La2} via a trace function on the generalized Hecke algebra of type B, which recovers HOM(ST), and the other is of the Kauffman bracket type, constructed in \cite{D0}, and which recovers KBSM(ST).

\subsection{Topological and algebraic set up}\label{tsetup}

\begin{figure}[H]
\begin{center}\includegraphics[width=3.1in]{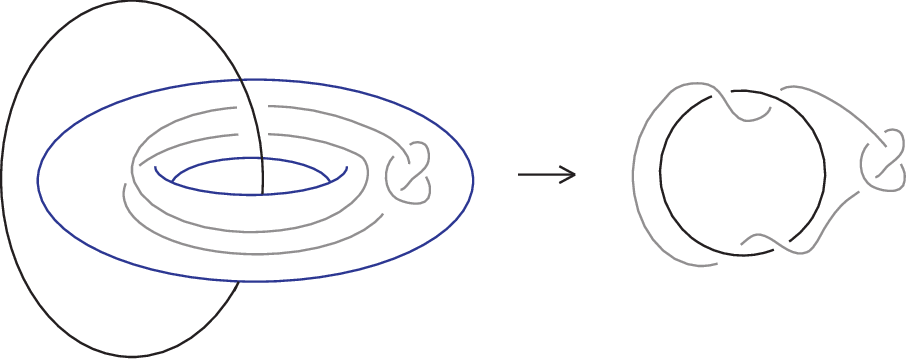}
\end{center}
\caption{A link in ST as a mixed link in $S^3$.}
\label{mli}
\end{figure}

We consider ST to be the complement of a solid torus in $S^3$. Then, an oriented link $L$ in ST can be represented by an oriented \textit{mixed link} in $S^{3}$, that is, a link in $S^{3}$ consisting of the unknotted fixed part $\widehat{I}$ representing the complementary solid torus in $S^3$, and the moving part $L$ that links with $\widehat{I}$. A \textit{mixed link diagram} is a diagram $\widehat{I}\cup \widetilde{L}$ of $\widehat{I}\cup L$ on the plane of $\widehat{I}$, where this plane is equipped with the top-to-bottom direction of $I$ (for an illustration see Figure~\ref{mli}).

\begin{figure}[H]
\begin{center}\includegraphics[width=2.7in]{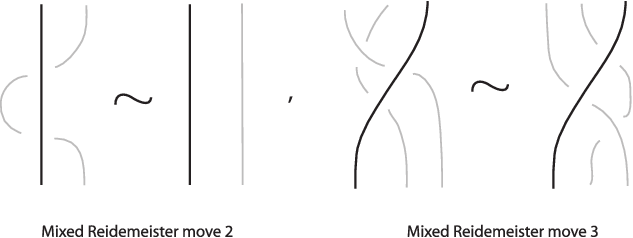}
\end{center}
\caption{The mixed Reidemeister moves.}
\label{mr}
\end{figure}

Isotopy in ST is then translated to mixed link isotopy in $S^3$ as follows:

\begin{theorem}[Theorem 5.2 \cite{LR1}]
Two links in ST are isotopic if and only if any two corresponding mixed link diagrams of theirs differ by planar isotopy and a finite sequence of the Reidemeister moves for the standard part of the mixed link and the mixed Reidemeister moves that involve both the fixed and the moving part of the mixed links and which are illustrated in Figure~\ref{mr}.
\end{theorem}

We now recall the analogue of the Alexander theorem for knots in ST (cf. Thm.~1 \cite{La1}):

\begin{theorem}[{\bf The analogue of the Alexander theorem for ST}]
A mixed link diagram $\widehat{I}\cup \widetilde{L}$ of $\widehat{I}\cup L$ may be turned into a \textit{mixed braid} $I\cup \beta$ with isotopic closure.
\end{theorem}

\begin{figure}[H]
\begin{center}\includegraphics[width=2.4in]{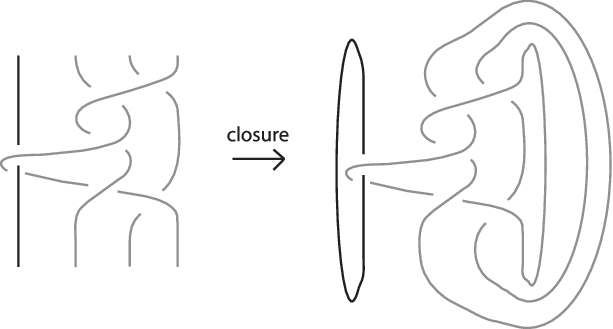}
\end{center}
\caption{The closure of a mixed braid to a mixed link.}
\label{cl}
\end{figure}

By the Alexander theorem for knots in solid torus, a mixed link diagram $\widehat{I}\cup \widetilde{L}$ of $\widehat{I}\cup L$ may be turned into a \textit{mixed braid} $I\cup \beta $ with isotopic closure (for an illustration see Figure~\ref{cl}). This is a braid in $S^{3}$ where, without loss of generality, its first strand represents $\widehat{I}$, the fixed part, and the other strands, $\beta$, represent the moving part $L$. The subbraid $\beta$ shall be called the \textit{moving part} of $I\cup \beta $.

\begin{figure}[H]
\begin{center}
\includegraphics[width=4.1in]{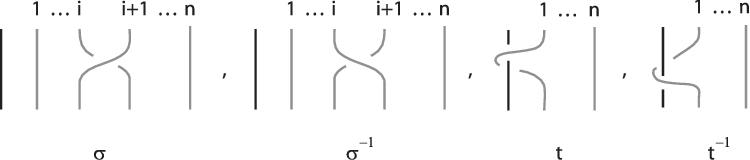}
\end{center}
\caption{The generators of $B_{1,n}$.}
\label{gen}
\end{figure}

The sets of braids related to ST form groups, which are in fact the Artin braid groups type B, denoted $B_{1,n}$, with presentation:

\[ B_{1,n} = \left< \begin{array}{ll}  \begin{array}{l} t, \sigma_{1}, \ldots ,\sigma_{n-1}  \\ \end{array} & \left| \begin{array}{l}
\sigma_{1}t\sigma_{1}t=t\sigma_{1}t\sigma_{1} \ \   \\
 t\sigma_{i}=\sigma_{i}t, \quad{i>1}  \\
{\sigma_i}\sigma_{i+1}{\sigma_i}=\sigma_{i+1}{\sigma_i}\sigma_{i+1}, \quad{ 1 \leq i \leq n-2}   \\
 {\sigma_i}{\sigma_j}={\sigma_j}{\sigma_i}, \quad{|i-j|>1}  \\
\end{array} \right.  \end{array} \right>, \]

\noindent where the generators $\sigma _{i}$ and $t$ are illustrated in Figure~\ref{gen}.

Let now $\mathcal{L}$ denote the set of oriented knots and links in ST. Isotopy in ST is then translated on the level of mixed braids by means of the following theorem:

\begin{theorem}[Theorem~4, \cite{LR2}] \label{markov}
 Let $L_{1} ,L_{2}$ be two oriented links in ST and let $I\cup \beta_{1} ,{\rm \; }I\cup \beta_{2}$ be two corresponding mixed braids in $S^{3}$. Then $L_{1}$ is isotopic to $L_{2}$ in ST if and only if $I\cup \beta_{1}$ is equivalent to $I\cup \beta_{2}$ in $\mathcal{B}$ by the following moves:
\[ \begin{array}{clll}
(i)  & Conjugation:         & \alpha \sim \beta^{-1}\, \alpha\, \beta, & {\rm if}\ \alpha ,\beta \in B_{1,n}. \\
(ii) & Stabilization\ moves: &  \alpha \sim \alpha\, \sigma_{n}^{\pm 1} \in B_{1,n+1}, & {\rm if}\ \alpha \in B_{1,n}. \\
(iii) & Loop\ conjugation: & \alpha \sim t^{\pm 1}\, \alpha\, t^{\mp 1}, & {\rm if}\ \alpha \in B_{1,n}. \\
\end{array}
\]

\end{theorem}

For an illustration of the moves in Theorem~\ref{markov}, see Figure~\ref{mbeq}.

\subsection{The HOMFLYPT skein module of ST}\label{homstpoly}

In \cite{HK} the HOMFLYPT skein module of the solid torus has been computed using diagrammatic methods by means of the following theorem:

\begin{theorem}[Kidwell--Hoste] \label{turaev}
The HOMFLYPT skein module skein module of the solid torus, HOM(ST), is a free, infinitely generated $\mathbb{Z}[u^{\pm1},z^{\pm1}]$-module isomorphic to the symmetric tensor algebra $SR\widehat{\pi}^0$, where $\widehat{\pi}^0$ denotes the conjugacy classes of non trivial elements of $\pi_1(\rm ST)$.
\end{theorem}

\begin{figure}[!ht]
\begin{center}
\includegraphics[width=1.2in]{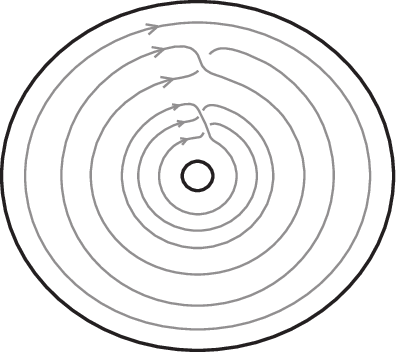}
\end{center}
\caption{A basic element of HOM(ST).}
\label{tur}
\end{figure}

A basic element of HOM(ST) in the context of \cite{HK}, is illustrated in Figure~\ref{tur}. Note that in the diagrammatic setting of \cite{HK}, ST is considered as ${\rm Annulus} \times {\rm Interval}$. 

HOM(ST) is well-studied and understood by now. It forms a commutative algebra with multiplication induced by embedding two solid tori in one in a standard way. Let now $\mathcal{B}^+$ denote the sub-algebra of HOM(ST), freely generated by elements that are clockwise oriented and let $\mathcal{B}^-$ denote the sub-algebra, freely generated by elements with counter-clockwise orientation. Let also $\mathcal{B}_k^+$ denote the sub-module generated by elements in $\mathcal{B}^+$ whose winding number is equal to $k\in \mathbb{N}$ and $\mathcal{B}_{-k}^-$ denote the sub-module generated by elements in $\mathcal{B}^-$ whose winding number is equal to $k$. As a linear space, $\mathcal{B}^+$ is graded by 
\[
\mathcal{B}^+ \cong \underset{k\geq 0}{\oplus}\, \mathcal{B}_k^+
\]

\noindent and similarly, $\mathcal{B}^-$ is graded by $$\mathcal{B}^- \cong \underset{k\geq 0}{\oplus}\, \mathcal{B}_{-k}^-.$$
Finally, it is worth mentioning the following module decomposition:
\[
{\rm HOM(ST)}\ =\ \underset{\lambda, \mu\geq 0}{\oplus}\, \mathcal{B}_{-\lambda}^-\, \otimes\, \mathcal{B}_{\mu}^+.
\]

\bigbreak

We now recover HOM(ST) using the braid approach following \cite{La2}. We consider the {\it generalized Hecke algebra} of type B, $\textrm{H}_{1, n}$ as the quotient of ${\mathbb C}\left[q^{\pm 1} \right]B_{1,n}$ over the quadratic relations $\sigma_{i}^2=(q-1)\sigma_{i}+q$. Namely:

\begin{equation*}
\textrm{H}_{1,n}(q)= \frac{{\mathbb C}\left[q^{\pm 1} \right]B_{1,n}}{ \langle \sigma_i^2 -\left(q-1\right)\sigma_i-q \rangle}.
\end{equation*}

\begin{figure}[H]
\begin{center}
\includegraphics[width=2.5in]{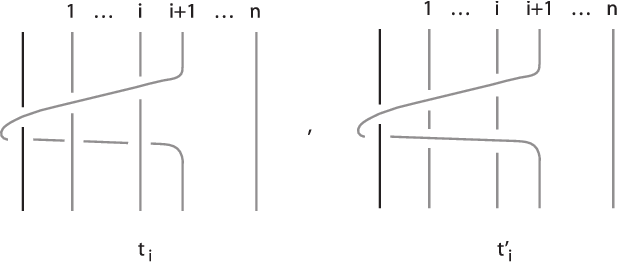}
\end{center}
\caption{The looping generators.}
\label{lp}
\end{figure}

This algebra is infinite dimensional, since the ``looping'' generator $t$ satisfies no polynomial relation. Note also that in \cite{La2} the algebra $\textrm{H}_{1,n}(q)$ is denoted as $\textrm{H}_{n}(q, \infty)$. From \cite{La2} we have:

\begin{theorem}{\cite[Proposition~1 \& Theorem~1]{La2}} \label{basesH}
The following sets form linear bases for ${\rm H}_{1,n}(q)$:
\[
\begin{array}{llll}
 (i) & \Sigma_{n} & = & \{t_{i_{1} } ^{k_{1} } \ldots t_{i_{r}}^{k_{r} } \cdot \sigma \} ,\ {\rm where}\ 0\le i_{1} <\ldots <i_{r} \le n-1,\\
 (ii) & \Sigma^{\prime} _{n} & = & \{ {t^{\prime}_{i_1}}^{k_{1}} \ldots {t^{\prime}_{i_r}}^{k_{r}} \cdot \sigma \} ,\ {\rm where}\ 0\le i_{1} < \ldots <i_{r} \le n-1, \\
\end{array}
\]
\noindent where $k_{1}, \ldots ,k_{r} \in {\mathbb Z}$ and $\sigma$ a basic element in ${\rm H}_{n}(q)$ and the ``looping generators'' $t_i$ and $t_i^{\prime}$ are defined as follows (for an illustration see Figure~\ref{lp}):

\begin{equation}
t_0^{\prime}\ =\ t_0\ :=\ t, \quad t_i^{\prime}\ =\ \sigma_i\ldots \sigma_1t\sigma_1^{-1}\ldots \sigma_i^{-1} \quad {\rm and}\quad t_i\ =\ \sigma_i\ldots \sigma_1t\sigma_1\ldots \sigma_i
\end{equation}
\end{theorem}

\smallbreak

We shall denote 
\begin{equation}
\Sigma\ :=\ \bigcup_n \Sigma_n \quad {\rm and\ similarly}\quad \Sigma^{\prime}\ :=\ \bigcup_n \Sigma_n^{\prime}.
\end{equation}

In \cite{La2} the bases $\Sigma^{\prime}_{n}$ are used for constructing a Markov trace on $\bigcup _{n=1}^{\infty }\textrm{H}_{1,n}(q)$.

\begin{theorem}{\cite[Theorem~6]{La2}} \label{tr}
Given $z,s_{k}$, with $k\in {\mathbb Z}$ specified elements in $R={\mathbb Z}\left[q^{\pm 1} \right]$, there exists
a unique linear Markov trace function
\begin{equation*}
{\rm tr}:\bigcup _{n=1}^{\infty }{\rm H}_{1,n}(q)  \to R\left(z,s_{k} \right),k\in {\mathbb Z}
\end{equation*}

\noindent determined by the rules:
\[
\begin{array}{lllll}
(1) & {\rm tr}(ab) & = & {\rm tr}(ba) & \quad {\rm for}\ a,b \in {\rm H}_{1,n}(q) \\
(2) & {\rm tr}(1) & = & 1 & \quad {\rm for\ all}\ {\rm H}_{1,n}(q) \\
(3) & {\rm tr}(ag_{n}) & = & z{\rm tr}(a) & \quad {\rm for}\ a \in {\rm H}_{1,n}(q) \\
(4) & {\rm tr}(a{t^{\prime}_{n}}^{k}) & = & s_{k}{\rm tr}(a) & \quad {\rm for}\ a \in {\rm H}_{1,n}(q),\ k \in {\mathbb Z}. \\
\end{array}
\]
\end{theorem}

Note that the use of the looping elements $t_i^{\prime}$ enable the trace ${\rm tr}$ to be defined by just extending  by rule (4) the three rules of the Ocneanu trace on the algebras ${\rm H}_n(q)$ \cite{Jo}. Using $\textrm{tr}$ a universal HOMFLYPT-type invariant for oriented links in ST is constructed. Namely, let $\mathcal{L}$ denote the set of oriented links in ST. Then:

\begin{theorem}{\cite[Definition~1]{La2}} \label{inv}
The function $X:\mathcal{L}$ $\rightarrow R(z,s_{k})$

\begin{equation*}
X_{\widehat{\alpha}} = \Delta^{n-1}
{\rm tr}\left(\pi \left(\alpha \right)\right),
\end{equation*}

\noindent where $\Delta:=-\frac{1-\lambda q}{\sqrt{\lambda } \left(1-q\right)} \left(\sqrt{\lambda } \right)^{e}$, $\lambda := \frac{z+1-q}{qz}$, $\alpha \in B_{1,n}$ is a word in the $\sigma _{i}$'s and $t^{\prime}_{i} $'s, $\widehat{\alpha}$ is the closure of $\alpha$, $e$ is the exponent sum of the $\sigma _{i}$'s in $\alpha $, and $\pi$ the canonical map of $B_{1,n}$ on ${\rm H}_{1,n}(q)$, such that $t\mapsto t$ and $\sigma _{i} \mapsto g_{i} $, is an invariant of oriented links in {\rm ST}.
\end{theorem}

\begin{remark}\rm
As shown in \cite{La2, DL2} the invariant $X$ recovers the HOMFLYPT skein module of ST. For a survey on the HOMFLYPT skein module of the lens spaces $L(p, 1)$ via braids, the reader is referred to \cite{DL3, DGLM}.
\end{remark}

In the braid setting of \cite{La2}, the elements of HOM(ST) correspond bijectively to the elements of the following set $\Lambda^{\prime}$:

\begin{equation}\label{Lpr}
\Lambda^{\prime}=\{ {t^{k_0}}{t^{\prime}_1}^{k_1} \ldots
{t^{\prime}_n}^{k_n}, \ k_i \in \mathbb{Z}\setminus\{0\}, \ k_i \geq k_{i+1}\ \forall i,\ n\in \mathbb{N} \}.
\end{equation}

\begin{figure}[H]
\begin{center}
\includegraphics[width=0.5in]{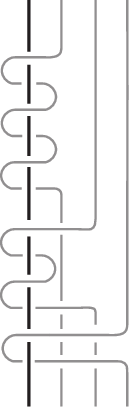}
\end{center}
\caption{The basic element $t^3{t_1^{\prime}}^2t_2^{\prime}$ of $\Lambda^{\prime}$.}
\label{belt}
\end{figure}

\noindent So, we have that $\Lambda^{\prime}$ is a basis of HOM(ST) in terms of braids (for an illustration see Figure~\ref{belt}). Note that $\Lambda^{\prime}$ is a subset of $\bigcup_n{\textrm{H}_{1,n}}$ and, in particular, $\Lambda^{\prime}$ is a subset of $\Sigma^{\prime}$. Note also that in contrast to elements in $\Sigma^{\prime}$, the elements in $\Lambda^{\prime}$ have no gaps in the indices, the exponents are ordered and there are no `braiding tails'. 

\begin{remark}\rm
It is worth mentioning that the invariant $X$ satisfies the following skein relation:
\[
\frac{1}{\sqrt{q}\sqrt{\lambda}}X_{L_+}-\sqrt{q}\sqrt{\lambda}X_{L_-}\ =\ \left(\sqrt{q}-\frac{1}{\sqrt{q}} \right) X_{L_o}.
\]
\end{remark}

\subsection{The Kauffman bracket skein module of ST}\label{kbsmst}

\begin{figure}[H]
\begin{center}
\includegraphics[width=4.2in]{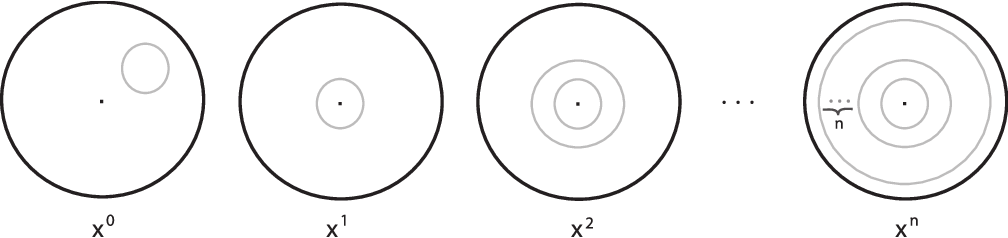}
\end{center}
\caption{The standard basis of KBSM(ST).}
\label{tur1}
\end{figure}

In this subsection we present results on the Kauffman bracket skein module of ST, KBSM(ST). A basis for KBSM(ST) is presented in \cite{HK} using diagrammatic methods (see also \cite{Tu}). More precisely:

\begin{theorem}[\cite{HK}]
The Kauffman bracket skein module of ST, KBSM(ST), is freely generated by an infinite set of generators $\left\{x^n\right\}_{n=0}^{\infty}$, where $x^n$ denotes $n$ parallel copies of a longitude of ST and $x^0$ is the affine unknot (see Figure~\ref{tur1}).
\end{theorem}

\begin{figure}[H]
\begin{center}
\includegraphics[width=0.7in]{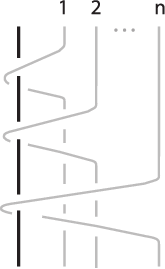}
\end{center}
\caption{The element $tt_1^{\prime}\ldots t_{n-1}^{\prime}$ of $B^{\prime}_{{\rm ST}}$.}
\label{ab1}
\end{figure}

In the braid setting, the elements of the standard basis of KBSM(ST) correspond bijectively to the elements of the following set (for an illustration see Figure~\ref{ab1}):

\begin{equation}\label{Lpr1}
{B}^{\prime}_{{\rm ST}}=\{ t{t^{\prime}_1} \ldots {t^{\prime}_n}, \ n\in \mathbb{N} \}.
\end{equation}

\begin{remark}\rm
Note now that the basis ${B}^{\prime}_{{\rm ST}}$ is not ``natural'', since on the level of braids, elements in ${B}^{\prime}_{{\rm ST}}$ have crossings. This was our main motivation for establishing a different and more suitable basis for KBSM(ST). We present this new basis in \S~\ref{nbas}. As we will see, this new basis is appropriate in extending the (universal) bracket polynomial, that recovers KBSM(ST), to an invariant of knots and links in the lens spaces $L(p,q)$.
\end{remark}

We now compute the Kauffman bracket skein module of ST via braids. Our starting point is the unique Markov trace constructed by Lambropoulou in \cite{La2} on the generalized Hecke algebra of type B, ${\rm H}_{1,n}$, and the universal invariant for knots and links in ST of the HOMFLYPT type (recall Theorems~\ref{tr} and \ref{inv}). 

\smallbreak

Following the work of V.F.R. Jones, we define the {\it generalized Temperley-Lieb algebra of type B} as a quotient of ${\rm H}_{1,n}$ over the ideal generated by $\sigma_{i, i+1},\, i\in \mathbb{N}\backslash \{0\}$, where
\begin{equation}\label{ideal}
\begin{array}{lcl}
\sigma_{i, i+1} & := & 1 + u\ (\sigma_i+ \sigma_{i+1}) + u^2\ (\sigma_i\sigma_{i+1}+\sigma_{i+1}\sigma_i)+u^3\ \sigma_i\sigma_{i+1}\sigma_i.\\
\end{array}
\end{equation}

More precisely, we have the following definition:

\begin{definition}\rm
The {\it generalized Temperley-Lieb algebra of type B}, $TL_{1, n}$, is defined as the quotient of the generalized Hecke algebra of type B, $H_{1,n}(q)$, over the ideal generated by the elements $\sigma_{i, i+1},\, i\in \mathbb{N}\backslash \{0\}$ defined in Eq.~\ref{ideal}.
\end{definition}

\begin{remark}\rm
Note that in \cite{FG} a different presentation for ${\rm H}_{1, n}$ is used, that involves the parameters $u, v$ and the quadratic relations 
\begin{equation}\label{quad}
{\sigma_{i}^2=(u-u^{-1})\sigma_{i}+1}.
\end{equation}

\noindent One can switch from one presentation to the other by a taking $\sigma_i = u \sigma_i$, $t = v t$ and $q = u^2$. We will adapt this setting from now on when working on $TL_{1, n}$.
\end{remark}

The necessary and sufficient conditions so that the Markov trace defined in ${\rm H}_{1, n}$ factors through to ${TL}_{1, n}$ are presented in \cite{D0}:

\begin{theorem}\label{tr2}
The trace defined in ${\rm H}_{1,n}$ factors through to $TL_{1, n}$ if and only if the trace parameters take the following values:
\begin{equation}\label{parameters}
z=-\frac{1}{u(1+u^2)}.
\end{equation}
\end{theorem} 

\begin{remark}\rm
Note now that for $z=-\, \frac{1}{u(1+u^2)}$, one deduces $\lambda\ =\ u^4$.
\end{remark}

We are now in position to present the universal invariant for knots and links in ST of the Kauffman bracket type \cite{D0}:

\begin{theorem}\label{inva}
The following invariant is the universal invariant of the Kauffman bracket type, for knots and links in ST:

\begin{equation*}
V_{\widehat{\alpha}}(u ,v)\ :=\ \left(-\frac{1+u^2}{u} \right)^{n-1} \left(u\right)^{2e} {\rm tr}\left(\overline{\pi} \left(\alpha \right)\right),
\end{equation*}

\noindent where $\alpha \in B_{1,n}$ is a word in the $\sigma _{i}$'s and $t^{\prime}_{i} $'s, $\widehat{\alpha}$ is the closure of $\alpha$, $e$ is the exponent sum of the $\sigma _{i}$'s in $\alpha $, $\overline{\pi}$ the canonical map of $B_{1,n}$ to ${\rm TL}_{1, n}$, such that $t\mapsto t$ and $\sigma _{i} \mapsto g_{i}$.
\end{theorem}

\begin{remark}\rm
The invariant $V$ recovers the Kauffman bracket skein module of ST since it gives distinct values to distinct elements of the basis of KBSM(ST), ${B}^{\prime}_{{\rm ST}}$. Indeed, $tr(tt_1^{\prime}\ldots t_n^{\prime})\, =\, s_1^n$.
\end{remark}

\begin{remark}\rm
It is worth emphasizing on the difference between the basis of the HOMFLYPT skein module of ST (recall Eq.~(\ref{Lpr})) and that of the Kauffman bracket skein module of ST (Eq.~(\ref{Lpr1})). For the Kauffman bracket there are no crossings, and this is due the skein relation in the definition of the Kauffman bracket. On the other hand, the skein relation that characterizes the HOMFLYPT polynomial does not annihilate the crossings, but changes its ``handiness'' instead.
\end{remark}

\begin{remark}\rm
\begin{itemize}
\item[i.] In \cite{D8} we introduce and study {\it knotoids} on the surface of the torus. Knotoids are open knotted curves in oriented surfaces $\Sigma$, that is, generic immersions of the unit interval $[0, 1]$ into $\Sigma$. In \cite{D8} we extend the notion of skein modules for knotoids and in particular, we show that the Kauffman bracket skein module of the Solid Torus is freely generated by elements of the form illustrated in Figure~\ref{difbas1}.

\begin{figure}[H]
\begin{center}
\includegraphics[width=3.8in]{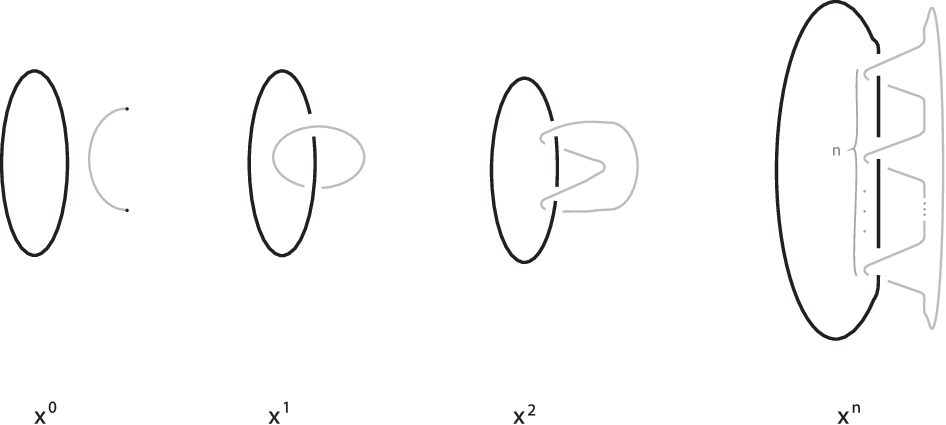}
\end{center}
\caption{A basis for KBSM(ST) for knotoids.}
\label{difbas1}
\end{figure}
\item[ii.] Finally, it is worth mentioning that in \cite{D7, D9} we construct the Kauffman bracket for pseudo knots in the Solid Torus and in the handlebody of genus $g$. Our intention is to generalize these invariants for other 3-manifolds, i.e., compute the corresponding skein module.
\end{itemize}
\end{remark}

\section{Skein modules of the lens spaces $L(p,1)$}\label{smlp}

\subsection{Topological and algebraic set up}\label{setup1}

\begin{figure}[H]
\begin{center}\includegraphics[width=3.3in]{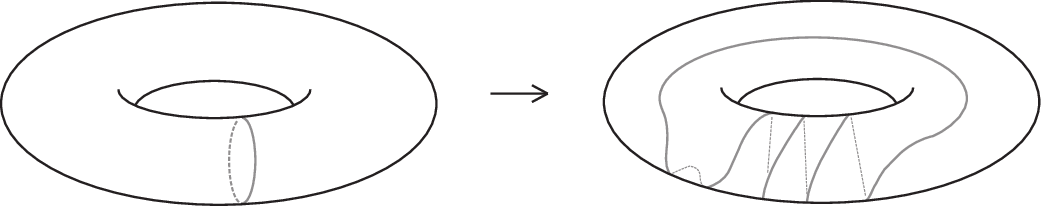}
\end{center}
\caption{The homeomorphism $h$.}
\label{st1}
\end{figure}

In this subsection we recall results from \cite{LR1, LR2, DL1}. We consider the lens spaces $L(p,1)$ obtained from $S^3$ by integral surgery along the unknot with coefficient $p$. Namely, we start from $S^3$ and we drill out a tubular neighborhood of the unknot, that is, a solid torus ${\rm ST}_1$. What is left is another solid torus ${\rm ST}_2$. We then glue these two solid tori via a homeomorphism $h$ on their boundaries, such that a meridian on $\partial {\rm ST}_1$ is mapped to a $(p, 1)$-curve on $\partial {\rm ST}_2$ (for an illustration see Figure~\ref{st1}).

\begin{remark}\rm
Note that a meridian bounds a disc in ST, and since this meridian is glued to a $(p, 1)$-torus knot, it follows that this $(p, 1)$-torus knot bounds a disc in $L(p,1)$.
\end{remark}

\begin{figure}[H]
\begin{center}\includegraphics[width=4.9in]{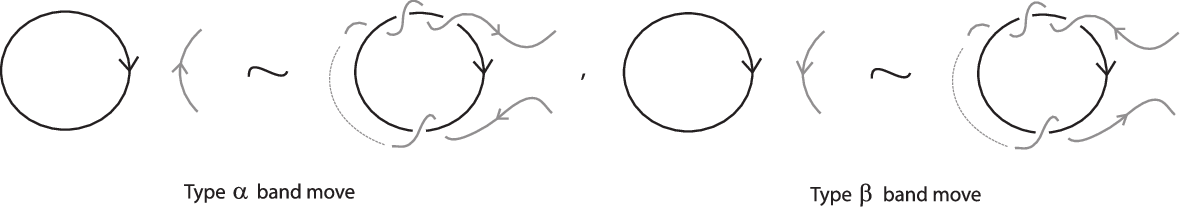}
\end{center}
\caption{The two types of band moves.}
\label{bm}
\end{figure}

Isotopy in $L(p, 1)$ can be viewed as isotopy in ST together with the {\it band moves} in $S^3$, which reflect the surgery description of the manifold, see Figure~\ref{bm} (see \cite{LR1}). 

In \cite{DL1} it is shown that it suffices to consider only one type of band move in order to describe isotopy for knots and links in
$L(p,1)$, and thus, isotopy between oriented links in $L(p, 1)$ is reflected in $S^3$ by means of the following theorem (cf.
\cite[Theorem 5.8]{LR1}, \cite[Theorem 6]{DL1}):

\begin{theorem}
Two oriented links in $L(p, 1)$ are isotopic if and only if two corresponding mixed link diagrams
of theirs differ by isotopy in ST together with a finite sequence of the type $\alpha$ (or $\beta$) band moves.
\end{theorem}

\begin{figure}[H]
\begin{center}
\includegraphics[width=5.1in]{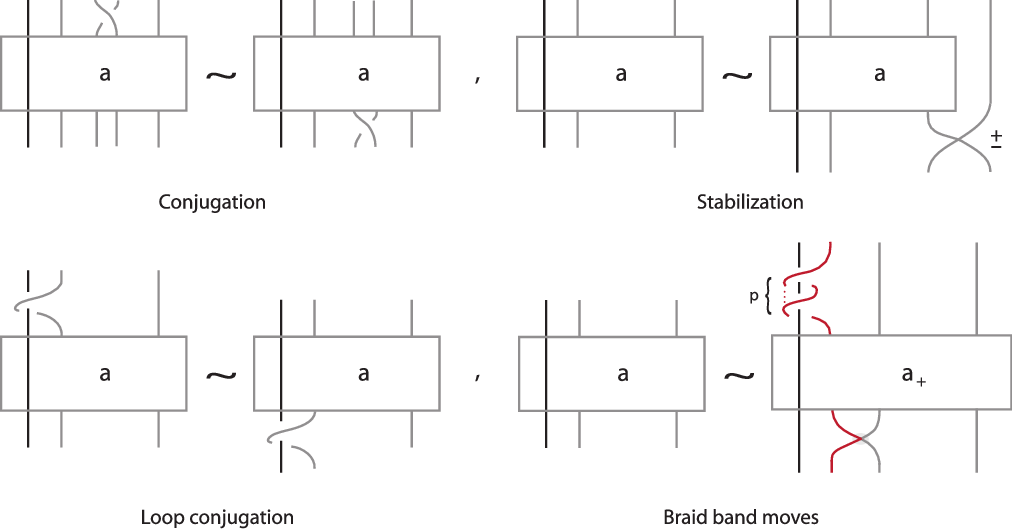}
\end{center}
\caption{Mixed braid equivalence.}
\label{mbeq}
\end{figure}

\begin{definition}\label{bbmdef}
We define a {\it braid band move} to be a move between mixed braids, which is a $\alpha$-type band move between their closures. It starts with a little band oriented downward, which, before sliding along a surgery strand, gets one twist {\it positive\/} or {\it negative\/} (see Figure~\ref{mbeq}). 
\end{definition}

\begin{remark}\rm
For now we shall only consider the $\alpha$-type band moves, since the result of an $\alpha$-type band move remain braided, while the result of a $\beta$-band move does not. Note also that there are two different types of braid band moves; the positive and the negative braid band move, depending on the kind of the twist the moving strand gets before sliding along the fixed strand. We denote a braid band move by bbm and, specifically, the result of a positive or negative braid band move performed on a mixed braid $\beta$ by $bbm_{\pm}(\beta)$.
\end{remark}

Isotopy in $L(p, 1)$ is then translated on the level of mixed braids by means of the following theorem:

\begin{theorem}[Theorem~5, \cite{LR2}] \label{markov1}
 Let $L_{1} ,L_{2}$ be two oriented links in $L(p,1)$ and let $I\cup \beta_{1} ,{\rm \; }I\cup \beta_{2}$ be two corresponding mixed braids in $S^{3}$. Then $L_{1}$ is isotopic to $L_{2}$ in $L(p,1)$ if and only if $I\cup \beta_{1}$ is equivalent to $I\cup \beta_{2}$ in $\mathcal{B}$ by the following moves:
\[ \begin{array}{clll}
(i)  & Conjugation:         & \alpha \sim \beta^{-1}\, \alpha\, \beta, & {\rm if}\ \alpha ,\beta \in B_{1,n}. \\
(ii) & Stabilization\ moves: &  \alpha \sim \alpha\, \sigma_{n}^{\pm 1} \in B_{1,n+1}, & {\rm if}\ \alpha \in B_{1,n}. \\
(iii) & Loop\ conjugation: & \alpha \sim t^{\pm 1}\, \alpha\, t^{\mp 1}, & {\rm if}\ \alpha \in B_{1,n}. \\
(iv) & Braid\ band\ moves: & \alpha \sim {t}^p\, \alpha_+\, \sigma_1^{\pm 1}, & a_+\in B_{1, n+1},\\
\end{array}
\]

\noindent where $\alpha_+$ is the word $\alpha$ with all indices shifted by +1 and $\beta_+$ is the word $\beta$ with all indices shifted by $q$. Note that moves (i), (ii) and (iii) correspond to link isotopy in {\rm ST} (recall Theorem~\ref{markov}).
\end{theorem}

\begin{figure}[H]
\begin{center}
\includegraphics[width=1.8in]{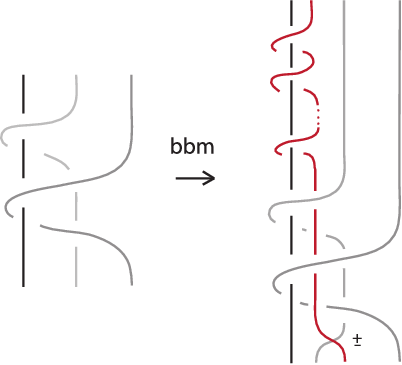}
\end{center}
\caption{Braid band move on the element $tt_1^{\prime}$.}
\label{bbmt}
\end{figure}

As noted before, our aim is to extend the universal invariants $X$ and $V$ for oriented (respectively framed) knots and links in ST, to invariants for knots and links in $L(p,1)$. Namely, we need to study the effect of bbm's to elements in a basis of the corresponding skein module of ST. Note now that, as shown in Figure~\ref{bbmt}, bbm's are not naturally described by elements in the bases $\Lambda^{\prime}$ for the HOM(ST), and by elements in ${B}^{\prime}_{{\rm ST}}$ for the case of the KBSM(ST). For this reason we first present more appropriate bases for both HOM(ST) and KBSM(ST).

\subsection{A new basis for HOM(ST)}\label{newbasishomst}

A different basis, $\Lambda$, for HOM(ST) is presented in \cite{DL2}. This basis is important for computing HOM($L(p,1)$), since, as noted before, it naturally describes the braid band moves, which extend the link isotopy in ST to link isotopy in $L(p,1)$. This was in fact our motivation for establishing this new basis. We have the following theorem:

\begin{theorem}{\cite[Theorem~2]{DL2}}\label{newbasis}
The following set forms a basis for HOM(ST):
\begin{equation}\label{basis}
\Lambda=\{t^{k_0}t_1^{k_1}\ldots t_n^{k_n},\ k_i \in \mathbb{Z}\setminus\{0\},\ k_i \geq k_{i+1}\ \forall i,\ n \in \mathbb{N} \}.
\end{equation}
\end{theorem}

\begin{figure}[!ht]
\begin{center}
\includegraphics[width=0.6in]{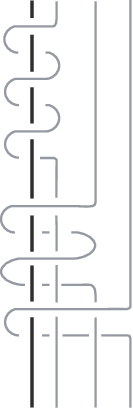}
\end{center}
\caption{The element $t^3t_1^2t_2$ in the new basis of HOM(ST).}
\label{basel}
\end{figure}

For an illustration of elements in the basis $\Lambda$ see Figure~\ref{basel}. Note that comparing the set $\Lambda$ with the set $\Sigma=\bigcup_n\Sigma_n$, we observe that in $\Lambda$ there are no gaps in the indices of the $t_i$'s and the exponents are in decreasing order. Also, there are no `braiding tails' in the words in $\Lambda$.

Our method for proving Theorem~\ref{newbasis} is the following: We first define total orderings in the sets $\Lambda^{\prime}$ and $\Lambda$, and we relate the two ordered sets via a lower triangular infinite matrix with invertible elements on the diagonal. More precisely, we start with elements in the basic set $\Lambda^{\prime}$ and we convert them into sums of elements in $\Sigma$. These elements consist of two parts: arbitrary monomials in the $t_i$'s followed by `braiding tails' in the bases of the algebras ${\rm H}_n(q)$. Then, these elements are converted into elements in the set $\Lambda$ by managing the gaps in the indices, by ordering the exponents of the $t_i$'s and by eliminating the `braiding tails'. These procedures are not independent since if for example we manage the gaps in the indices of the looping generators $t_i$'s, `braiding tails' may occur and also the exponents of the $t_i$'s may alter. Similarly, when the `braiding tails' are eliminated, gaps in the indices of the $t_i$'s might occur. This is a long procedure that eventually stops and only elements in the set $\Lambda$ remain. Using the change of basis matrix we finally show that the set $\Lambda$ is linearly independent and hence, the set $\Lambda$ forms a basis for HOM(ST).

\smallbreak

In order to define the ordering relation, we introduce the notion of the {\it index} of a word $w$, denoted $ind(w)$.

\begin{definition}
Let $w \in \Lambda$ (respectively in $\Lambda^{\prime}$). Then, $ind(w)$ is defined to be the highest index of the $t_i$'s (resp. of the $t_i^{\prime}$'s) in $w$. If $w \in \Sigma$ (or in $\Sigma^{\prime}$), then $ind(w)$ is defined as above by ignoring possible gaps in the indices of the looping generators and by ignoring the braiding parts in the algebras $\textrm{H}_{n}(q)$. Finally, we define the index of a monomial in $\textrm{H}_{n}(q)$ to be equal to $0$.
\end{definition}

	\begin{example} \rm
	\[
	\begin{array}{llllll}
{i.} & ind({t^{\prime}}^{k_0}{t^{\prime}_1}^{k_1}\ldots {{t^{\prime}}_m}^{k_m}) & = & m & = & ind(t^{k_0}{t_1}^{k_1}\ldots {t_m}^{k_m})\\
{ii.} & ind(t^{k_0}t_3^{k_3}\cdot \sigma)                                      & = & 3 & &,\ \sigma \in {\rm H}_n(q)\\
{iii.} & ind(\sigma)                                                           & = & 0 & &,\  \sigma \in {\rm H}_n(q)
	\end{array}
	\]
	\end{example}

We are now in position to state the ordering relation.

\begin{definition}{\cite[Definition~2]{DL2}} \label{order}
\rm
Let $w={t^{\prime}_{i_1}}^{k_1}\ldots {t^{\prime}_{i_{\mu}}}^{k_{\mu}}\cdot \beta_1$ and $u={t^{\prime}_{j_1}}^{\lambda_1}\ldots {t^{\prime}_{j_{\nu}}}^{\lambda_{\nu}}\cdot \beta_2$ in $\Sigma^{\prime}$, where $k_t , \lambda_s \in \mathbb{Z}$ for all $t,s$ and $\beta_1, \beta_2 \in H_n(q)$. Then, we define the following ordering in $\Sigma^{\prime}$:

\smallbreak

\begin{itemize}
\item[(a)] If $\sum_{i=0}^{\mu}k_i < \sum_{i=0}^{\nu}\lambda_i$, then $w<u$.

\vspace{.1in}

\item[(b)] If $\sum_{i=0}^{\mu}k_i = \sum_{i=0}^{\nu}\lambda_i$, then:

\vspace{.1in}

\noindent  (i) if $ind(w)<ind(u)$, then $w<u$,

\vspace{.1in}

\noindent  (ii) if $ind(w)=ind(u)$, then:

\vspace{.1in}

\noindent \ \ \ \ ($\alpha$) if $i_1=j_1, \ldots , i_{s-1}=j_{s-1}, i_{s}<j_{s}$, then $w>u$,

\vspace{.1in}

\noindent \ \ \  ($\beta$) if $i_t=j_t$ for all $t$ and $k_{\mu}=\lambda_{\mu}, k_{\mu-1}=\lambda_{\mu-1}, \ldots, k_{i+1}=\lambda_{i+1}, |k_i|<|\lambda_i|$, then $w<u$,

\vspace{.1in}

\noindent \ \ \  ($\gamma$) if $i_t=j_t$ for all $t$ and $k_{\mu}=\lambda_{\mu}, k_{\mu-1}=\lambda_{\mu-1}, \ldots, k_{i+1}=\lambda_{i+1}, |k_i|=|\lambda_i|$ and $k_i>\lambda_i$, then $w<u$,

\vspace{.1in}

\noindent \ \ \ \ ($\delta$) if $i_t=j_t\ \forall t$ and $k_i=\lambda_i$, $\forall i$, then $w=u$.

\end{itemize}

The ordering in the set $\Sigma$ is defined as in $\Sigma^{\prime}$, where $t_i^{\prime}$'s are replaced by $t_i$'s.
\end{definition}

\begin{example} \rm
\[
\begin{array}{rlllll}
{\rm i.} & {t}^3{t_1^{\prime}}^4 & < &  {t}^4{t_1^{\prime}}^{5} & {\rm since}\ 3+7\ <\ 4+5\ &,\ {\rm \left(Def.~\ref{order}a\right)}\\
&&&&&\\
{\rm ii.} & {t}{t_1^{\prime}}^2{t_2^{\prime}}^3 & < &  {t}^2{t_1^{\prime}}^2{t_2^{\prime}}{t_3^{\prime}} & {\rm since}\ ind({t}{t_1^{\prime}}^2{t_2^{\prime}}^3)<ind({t}^2{t_1^{\prime}}^2{t_2^{\prime}}{t_3^{\prime}}) &,\ {\rm \left(Def.~\ref{order}b(i)\right)}\\
&&&&&\\
{\rm iii.} & {t}^2{t_1^{\prime}}^3{t_3^{\prime}}^4{t_4^{\prime}}^8 & < &  {t}^8{t_1^{\prime}}{t_2^{\prime}}{t_4^{\prime}}^7 & {\rm since}\ 3\ =\ i_3\ >\ j_3\ =\ 2 &,\ {\rm \left(Def.~\ref{order}b(ii)(\alpha)\right)}\\
&&&&&\\
{\rm iv.} & {t}^2{t_1^{\prime}}^4{t_3^{\prime}}{t_4^{\prime}}^3 & < &   {t}^{14}{t_1^{\prime}}^{-8}{t_3^{\prime}}{t_4^{\prime}}^3 & {\rm since}\ |4|\ <\ |-8| &,\ {\rm \left(Def.~\ref{order}b(ii)(\beta)\right)}\\
&&&&&\\
{\rm v.} & t{t_1^{\prime}}^4{t_3^{\prime}}{t_4^{\prime}}^3 & < &  {t}^{10}{t_1^{\prime}}^{-4}{t_3^{\prime}}{t_4^{\prime}}^3 & {\rm since}\ -4\ <\ 4 &,\ {\rm \left(Def.~\ref{order}b(ii)(\gamma)\right)}
\end{array}
\]
\end{example}

\begin{remark}\label{nt} \rm
We shall denote a monomial of the form $t_i^{k_i}\ldots t^{k_{i+m}}_{i+m}$ by $\tau_{i,i+m}^{k_{i,i+m}}$, and similarly ${\tau^{\prime}}_{i,i+m}^{k_{i,i+m}}:={t^{\prime}}_i^{k_i}\ldots {t^{\prime}}^{k_{i+m}}_{i+m}$, for $m\in \mathbb{N}$, $k_j\neq 0$ for all $j$.
\end{remark}

In order to relate the two basic sets $\Lambda^{\prime}$ and $\Lambda$ via an infinite matrix, we need to define the \textit{subsets of level $k$}, $\Lambda_{(k)}$ and $\Lambda^{\prime}_{(k)}$, of $\Lambda$ and $\Lambda^{\prime}$ respectively ({\cite[Definition~3]{DL2}}), to be the sets:

\begin{equation}\label{levk}
\begin{array}{l}
\Lambda_{(k)}:=\{t_0^{k_0}t_1^{k_1}\ldots t_{m}^{k_m} | \sum_{i=0}^{m}{k_i}=k,\ k_i \in \mathbb{Z}\setminus\{0\},\  k_i \leq k_{i+1}\ \forall i \}\\
\\
\Lambda^{\prime}_{(k)}:=\{{t^{\prime}_0}^{k_0}{t^{\prime}_1}^{k_1}\ldots {t^{\prime}_{m}}^{k_m} | \sum_{i=0}^{m}{k_i}=k,\ k_i \in \mathbb{Z}\setminus\{0\},\  k_i \leq k_{i+1}\ \forall i \}
\end{array}
\end{equation}

\noindent In \cite{DL2} it is shown that the sets $\Lambda_{(k)}$ and $\Lambda^{\prime}_{(k)}$ are totally ordered and well ordered for all $k$ (\cite[Propositions~1 \& 2]{DL2}). Note that in \cite{DLP} the exponents in the monomials of $\Lambda$ are in decreasing order, while here the exponents are considered in increasing order, which is totally symmetric. 

\bigbreak

We finally define the set $\Lambda^{aug}$, which augments the basis $\Lambda$ and its subset of level $k$, and we also introduce the notion of \textit{homologous words}.

\begin{definition}\rm
We define the set:
\begin{equation}\label{lamaug}
\Lambda^{aug}\ :=\{t_0^{k_0}t_1^{k_1}\ldots t_{n}^{k_n},\ k_i \in \mathbb{Z}\backslash \{0\}\}.
\end{equation}
\noindent and the {\it subset of level} $k$, $\Lambda^{aug}_{(k)}$, of $\Lambda^{aug}$:

\begin{equation}\label{lamaug2}
\Lambda^{aug}_{(k)}:=\{t_0^{k_0}t_1^{k_1}\ldots t_{m}^{k_m} | \sum_{i=0}^{m}{k_i}=k,\ k_i \in \mathbb{Z}\backslash \{0\}\}
\end{equation}
\end{definition}

\begin{definition}\label{homw} \rm
We shall say that two words $w^{\prime}\in \Lambda^{\prime}$ and $w\in \Lambda$ are {\it homologous}, denoted $w^{\prime}\sim w$, if $w$ is
obtained from $w^{\prime}$ by turning $t^{\prime}_i$ into $t_i$ for all $i$.
\end{definition}

\begin{example}\rm
The words $t^{-2}{t_1^{\prime}}^{2}{t_2^{\prime}}$ and $t^{-2}{t_1}^{2}{t_2}$ are homologous. Note also that $\Lambda^{\prime} \ni t \ \sim\ t\in \Lambda$.
\end{example}
	
\bigbreak

In order to relate the sets $\Lambda^{\prime}$ and $\Lambda$ via a lower triangular matrix with invertible elements in the diagonal, we first express elements in $\Lambda^{\prime}$ to elements containing the $t_i$'s. Indeed we have that (Theorem 7, \cite{DL2}):

\begin{equation}\label{convert2}
\begin{array}{lcl}
t^{k_0}{t_1^{\prime}}^{k_1} \ldots {t_m^{\prime}}^{k_m} & = & q^{-\underset{n=1}{\overset{m}{\sum}}\, {nk_n}}\cdot\ t^{k_0}t_1^{k_1}\ldots t_m^{k_m} \ + \ \underset{i}{\sum}\, {f_i(q)\cdot t^{k_0}t_1^{k_1}\ldots t_m^{k_m}\cdot w_i} \ +\\
&&\\
& + & \underset{j}{\sum}\, {g_j(q)\tau_j \cdot u_j},
\end{array}
\end{equation}
\noindent where $w_i, u_j \in {\rm H}_{m+1}(q), \forall i$, $\tau_j \in \Sigma_n$, such that $\tau_j < t^{k_0}t_1^{k_1}\ldots t_m^{k_m}, \forall j$ and $f_i, g_j \in \mathbb{C}$, for all $i, j$.

\begin{example}\rm
We shall now give an example of a monomial in $\Lambda^{\prime}$ converted into sums of elements in $\Sigma$. We have that $t^2{t_1^{\prime}}\ = \ q^{-1}\, t^2{t_1}\ +\ (q^{-1}-1)\, t^3\, \sigma_1$. For more details the interested reader should refer to the technical lemmas \cite[Lemmas~3, 4, 5, 6, 9 \& 11]{DL2}.
\end{example}

Equations~(\ref{convert2}) suggest that an element in $\Lambda^{\prime}$ can be written as a sum of elements in $\Lambda$, where one term is the homologous word, another term is the homologous word followed by a braiding ``tail'', and all other terms in the sum consist of lower order terms followed by braiding ``tails''. These elements belong to $\Sigma_n$ since they may have gaps in their indices. We then manage the gaps in the indices (Theorem~8, \cite{DL2}), namely:

\begin{equation}\label{gaps}
\Sigma_n \ni \tau\ \widehat{=}\ \underset{i}{\sum}\, f_{i}(q)\, \tau_i\cdot w_i\ :\ \tau_i\in \Lambda^{aug},\ w_i\in {\rm H}_{n}(q),\ \forall i,
\end{equation}

\noindent where $\widehat{=}$ denotes that conjugation is applied in the process.

\smallbreak

Equation~(\ref{gaps}) is best demonstrated in the following example on a word with two gaps in the indices of the loop generators.

\begin{example}\rm For the 2-gap word $t^{k_0}t_2t_4\in \Sigma$ we have:
\[
\begin{array}{lclclc}
t^{k_0}\underline{t_2}t_4 & = & \underline{t^{k_0} \sigma_2}\, t_1\, \underline{\sigma_2\, t_4} & = & \sigma_2 t^{k_0}t_1t_4 \sigma_2 &\widehat{=} \\
&&&&&\\ 
& \widehat{=} & t^{k_0}t_1\underline{t_4} \sigma_2^2 & = & \underline{t^{k_0}t_1\, \sigma_4\sigma_3}\, t_2\, \sigma_3\sigma_4\, \sigma_2^2 & =\\
&&&&&\\ 
& = & \sigma_4\sigma_3\, t^{k_0}t_1 t_2\, \sigma_3\sigma_4\, \sigma_2^2 &  \widehat{=} & t^{k_0}t_1 t_2\, \sigma_3\sigma_4\, \sigma_2^2\, \sigma_4\sigma_3&
\end{array}
\]
\end{example}

We now deal with the elements in $\Lambda^{aug}_{(k)}$ that are followed by a braiding ``tail'' $w$ in $H_n(q)$. More precisely we have (Theorem~9, \cite{DL2}):

$$\tau \cdot w \ \widehat{\cong}\  \sum_{j}{f_j(q,z)\cdot \tau_j},$$

\noindent such that $\tau_j\in \Lambda^{aug}_{(k)}$ and $\tau_{j} < \tau$, for all $j$.

\begin{example}\rm In this example we demonstrate how to eliminate the `braiding tail' in a word.
\[
\begin{array}{lclclclclc}
t^{k}\underline{t_1^2t_2}\sigma_1^{-1} & = & t^{k}t_1t_2\underline{t_1\sigma_1^{-1}} & = & t^{k}t_1t_2\sigma_1\underline{t}& \widehat{=} & t^{k+1}\underline{t_1t_2}\sigma_1 & = & t^{k+1}t_2\underline{t_1\sigma_1} & =\\
&&&&&&&&&\\
& = & (q-1)t^{k+1}\underline{t_2t_1} & + & q t^{k+1}t_2 \sigma_1 \underline{t} & \widehat{=} & (q-1)t^{k+1}t_1t_2 & + & q t^{k+2}\underline{t_2} \sigma_1 & =\\
&&&&&&&&&\\
& = & (q-1)t^{k+1}t_1t_2 & + & q \underline{t^{k+2}\sigma_2}\, t_1\, \sigma_2 \sigma_1 & = & (q-1)t^{k+1}t_1t_2 & + & q \underline{\sigma_2} t^{k+2}\, t_1\, \sigma_2 \sigma_1 & \widehat{=}\\
&&&&&&&&&\\
& \widehat{=} & (q-1)t^{k+1}t_1t_2 & + & q t^{k+2}\, t_1\, \underline{\sigma_2 \sigma_1 \sigma_2} & = & (q-1)t^{k+1}t_1t_2 & + & q t^{k+2}\, t_1\, \underline{\sigma_2 \sigma_1 \sigma_2} & =\\
&&&&&&&&&\\
& = & (q-1)t^{k+1}t_1t_2 & + & q t^{k+2}\, t_1\, \sigma_1 \underline{\sigma_2} \sigma_1 & \cong & (q-1)t^{k+1}t_1t_2 & + & qz t^{k+2}\, t_1\, \underline{\sigma_1^2} & =\\
&&&&&&&&&\\
& = & (q-1)t^{k+1}t_1t_2 & + & q(q-1)z t^{k+2}\, \underline{t_1\, \sigma_1} & + & q^2z t^{k+2}\, t_1 & = & &\\
&&&&&&&&&\\
& = & (q-1)t^{k+1}t_1t_2 & + & q(q-1)^2z t^{k+2}\, t_1 & + &q^2(q-1)z t^{k+2}\, \underline{\sigma_1 t}  & + & q^2z t^{k+2}\, t_1 & \widehat{\cong}\\
&&&&&&&&&\\
& \widehat{\cong} & (q-1)t^{k+1}t_1t_2 & + & q(q-1)^2z t^{k+2}\, t_1 & + &q^2(q-1)z^2 t^{k+3}  & + & q^2z t^{k+2}\, t_1 & \widehat{\cong}\\
\end{array}
\]
\end{example}

\bigbreak

One very important result in \cite{DL2} is that one can change the order of the exponents by using conjugation and stabilization moves on elements in $\Lambda$ and express them as sums of monomials in $t_i$'s with arbitrary exponents and which are of lower order than the initial elements in $\Lambda$. Note that both conjugation and stabilization moves are captured by the trace rules, and that we translate here Theorem~10 in \cite{DL2} using the trace: 

$$tr(\tau_{0,m}^{k_{0,m}}\cdot w)\ =\ tr\left(\underset{j}{\sum}\, {\tau_{0,j}^{\lambda_{0,j}}\cdot w_j}\right),$$

\noindent where $\tau_{0,j}^{\lambda_{0,j}} \in \Lambda$ and $w, w_j \in \bigcup_{n\in \mathbb{N}}{\rm H}_n(q)$ for all $j$.

\begin{example}\rm 
Consider the element $tt_1^2t_2 \in \Lambda^{aug}$. We have that:
\[
\begin{array}{lclclclc}
t\underline{t_1^2}t_2 & = & t\, \sigma_1t\underline{\sigma_1^2} t\sigma_1\, t_2 & = & (q-1)\, t\, \underline{\sigma_1t\sigma_1} t\sigma_1\, t_2 & + & q\, \underline{t\, \sigma_1} t^2 \sigma_1\, t_2 & \widehat{=}\\
&&&&&&&\\
&& & \widehat{=} & (q-1)\, t\, t_1 t\sigma_1\, t_2 & + & q\, t^2 \underline{\sigma_1\, t_2} t\, \sigma_1 & =  \\
&&&&&&&\\
&& & = & (q-1)\, t^2\, t_1 \, t_2\, \sigma_1 & + & q\, t^2\, t_1 \, t_2 & 
\end{array}
\]
\end{example}

Finally, in \cite{DL2} it is shown that the infinite matrix converting elements of the basis $\Lambda^{\prime}$ of HOM(ST) to elements of the set $\Lambda$ is a block diagonal matrix, where each block corresponds to a subset of $\Lambda^{\prime}$ of level $k$ and it is an infinite lower triangular matrix with invertible elements in the diagonal. Using this infinite diagonal matrix, in \cite[Theorem~11]{DL2} it is shown that the set $\Lambda$ is linearly independent and hence, it forms a basis for HOM(ST).

\subsection{A new basis for KBSM(ST)}\label{nbas}

KBSM(ST) is particularly important since, as shown in \cite{P}, in order to compute KBSM($L(p,1)$), it suffices to consider elements in KBSM(ST) and study the effect of band moves on these elements. Namely,

\begin{equation}\label{przkmsm}
{\rm KBSM}\left(L(p, 1)\right)=\frac{{\rm KBSM}({\rm ST})}{<a-bbm(a)>}, \quad {\rm where}\ a\ {\rm basis\ element\ of\ KBSM(ST)}.
\end{equation}

\begin{figure}[H]
\begin{center}
\includegraphics[width=3.9in]{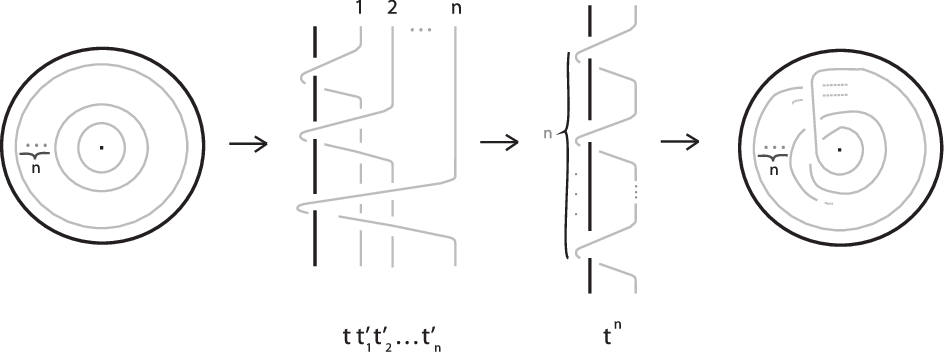}
\end{center}
\caption{Elements in the bases $B^{\prime}_{{\rm ST}}$ and $B_{\rm ST}$ of KBSM(ST).}
\label{ab}
\end{figure}

Equivalently, in order to compute KBSM$(L(p,1))$ we need to solve the infinite system of equations obtained by imposing on the universal invariant $V$, of the Kauffman bracket type, relations coming from the performance of braid band moves on elements in a basis of KBSM(ST). Toward that end, we first present a more suitable basis that involves elements with no crossings on the level of braids, and that it naturally describes bbm's (for an illustration see Figure~\ref{ab}). More precisely, we have the following result:

\begin{theorem}\label{newbasis1}
The following set is a basis for KBSM(ST):
\begin{equation}\label{basis1}
B_{\rm ST}\ =\ \{t^{n},\ n \in \mathbb{N} \}.
\end{equation}
\end{theorem}

In order to prove Theorem~\ref{newbasis1} in \cite{D0} we follow the following procedure (see also \cite{D1}):

\bigbreak

We first prove that the sets $B^{\prime}_{{\rm ST}}$ and $B_{{\rm ST}}$ are totally ordered and well-ordered sets when equipped with the ordering relation of Definition~\ref{order}.

\smallbreak

We then pass from the set $B^{\prime}_{{\rm ST}}$ consisting of monomials of the form $tt_1^{\prime}t_2^{\prime}\ldots t_n^{\prime},\ n\in \mathbb{N}$, to an augmented set ${B_{{\rm ST}}}^{aug}$, consisting of monomials of the form $t^n$, where $n\in \mathbb{Z}$. In particular in \cite{D0} it is shown that if $\tau\in B^{\prime}_{{\rm ST}}$ such that $in(\tau)=n$, then
\[
\tau \ \widehat{\underset{{\rm skein}}{\cong}} \ \underset{i\in \mathbb{Z}\, :\, |i|\, \leq\, n+1}{\sum}\, a_i\, t^i,
\]
\noindent where $a_i$ coefficients for all $i$.

\begin{example}\rm
In this example, we convert elements in the standard basis $B^{\prime}_{{\rm ST}}$ to sums of elements in ${B_{{\rm ST}}}^{aug}$.
\smallbreak

\qquad \qquad \qquad \qquad i. The case $tt_1^{\prime}$ is illustrated in Figure~\ref{bas1}. 

\begin{figure}[H]
\begin{center}
\includegraphics[width=5.1in]{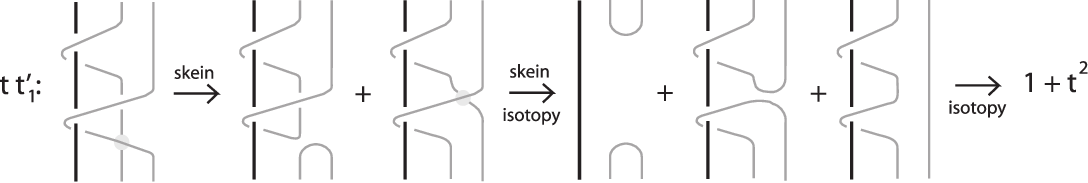}
\end{center}
\caption{Expressing $tt_1^{\prime}$ in terms of elements in ${B_{{\rm ST}}}^{aug}$.}
\label{bas1}
\end{figure}
\[
\begin{array}{llcl}
{\rm ii.} & t\, t_1^{\prime}\, t_2^{\prime} & \widehat{\underset{{\rm skein}}{\cong}} & t^{-1}\, +\, t\, +\, t^3\\
&&&\\
{\rm iii.} & t\, t_1^{\prime}\, t_2^{\prime}\, t_3^{\prime}\, t_4^{\prime}\, t_5^{\prime} & \widehat{\underset{{\rm skein}}{\cong}} & t^{-4}\, +\, t^{-2}\, +\, 1\, +\, t^2\, +\, t^4\, +\, t^6\\
&&&\\
{\rm iv.} & t\, t_1^{\prime}\, t_2^{\prime}\, t_3^{\prime}\, t_4^{\prime}\, t_5^{\prime}\, t_6^{\prime} & \widehat{\underset{{\rm skein}}{\cong}} & t^{-5}\, +\, t^{-3}\, +\, t^{-1}\, +\, t\, +\, t^3\, +\, t^5\, +\, t^7\\
\end{array}
\]
\end{example}

\smallbreak
We now express the elements $t^i \in {B_{{\rm ST}}}^{aug} \backslash B_{{\rm ST}}$ to sums of elements in $B_{{\rm ST}}$ (Proposition~2, \cite{D0}). In particular, we show that if $n\in \mathbb{N}$, then:
\[
t^{-n}\ \widehat{\underset{{\rm skein}}{\cong}}\ \underset{i=0}{\overset{n}{\sum}}\, a_i\, t^{i}.
\]

\begin{figure}[H]
\begin{center}
\includegraphics[width=4.8in]{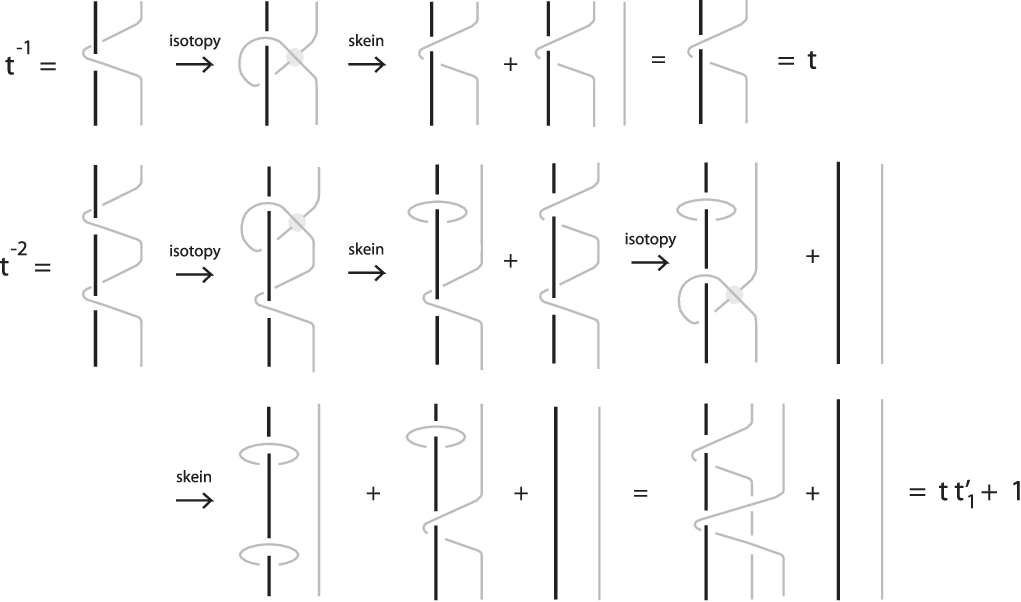}
\end{center}
\caption{Converting $t^{-1}$ and $t^{-2}$ to elements in $B^{\prime}_{{\rm ST}}$.}
\label{nexp}
\end{figure}

In Figure~\ref{nexp} we illustrate some examples of converting elements in ${B_{{\rm ST}}}^{aug} \backslash B_{{\rm ST}}$ to sums of elements in $B_{{\rm ST}}^{\prime}$.

The two (ordered) sets $B^{\prime}_{\rm ST}$ and $B_{\rm ST}$ are finally related via a lower triangular infinite matrix with invertible elements on the diagonal, and we conclude that the set $B_{\rm ST}$ forms a basis for KBSM(ST).

\subsection{Toward the HOMFLYPT skein module of the lens spaces $L(p,1)$}\label{homlp1}

In this subsection we go back to our original goal, which is the computation of the HOMFLYPT skein module of the lens spaces $L(p,1)$. As explained before, in order to compute HOM$(L(p,1))$, we must normalize the universal invariant $X$ of the HOMFLYPT type, by making it satisfy every possible braid band move. In \cite{DLP} it is shown that the performance of a {\it bbm} on a mixed braid in $B_{1, n}$ reduces to performing {\it bbm}'s on elements in the canonical basis, $\Sigma_n^{\prime}$, of the algebra ${\rm H}_{1,n}(q)$ and, in fact, on their first moving strand. Namely, if $\alpha \in B_{1, n}$ and $\beta_i \in \Sigma_n^{\prime}$, for all $i$ (\cite[Proposition~1]{DLP}). We have that:
\[
X_{\widehat{\alpha}}\, =\, X_{\widehat{bbm{(\alpha)}}}\ \Leftrightarrow\ \underset{i}{\sum}\, \left(X_{\widehat{\beta_i}}\, =\, X_{\widehat{bbm_1{(\beta_i)}}}\right)
\]

\noindent Namely, in order to compute HOM($L(p,1))$, it suffices to consider the performance of braid band moves on the first strand of elements in the set $\Sigma^{\prime}$. This simplifies the infinite system of equations. We then simplify the derived equations further, by showing that the equations obtained by performing bbm's on the first moving strand of elements in $\Sigma_n^{\prime}$ are equivalent to the equations obtained by performing bbm's on the first moving strand of elements in $\Sigma_n$ (\cite[Proposition~2]{DLP}). That is, 
\[
X_{\widehat{s^{\prime}}}\, =\, X_{\widehat{bbm_1{(s^{\prime})}}}\ \Leftrightarrow\ \underset{i}{\sum}\, \left(X_{\widehat{s_i}}\, =\, X_{\widehat{bbm_1{(s_i)}}}\right),
\]
\noindent where $s^{\prime}\in \Sigma_n^{\prime}$ and $s_i \in \Sigma_n$ for all $i$.

\smallbreak

Recall now that elements in $\Sigma$ consist of two parts: a monomial in $t_i$'s with possible gaps in the indices and unordered exponents, followed by a `braiding tail' in the basis of ${\rm H}_n(q)$. This allows us to further reduce the computation to elements in the basis $\Lambda$ of HOM(ST). In order to obtain elements in the augmented ${\rm H}_{n}(q)$-module $\Lambda^{aug}$ (followed by `braiding tails'), we first handle the gaps in the indices of the looping generators of elements in $\Sigma$. It should be noted that the performance of a bbm  is now seen as occurring on any moving strand (\cite[Proposition~3]{DLP}). Namely:
\[
X_{\widehat{s}}\, =\, X_{\widehat{bbm_1{(s)}}}\ \Leftrightarrow\ \underset{i}{\sum}\, \left(X_{\widehat{h_i}}\, =\, X_{\widehat{bbm{(h_i)}}}\right),
\]
\noindent where $s\in \Sigma_n$ and $h_i$ in the ${\rm H}_{n}(q)$-module $\Lambda^{aug}$.

\smallbreak

In \cite{DLP}, it shown that equations derived from elements in the $\Lambda^{aug}$-module ${\rm H}_{n}(q)$ are equivalent to the equations obtained by performing bbm's on an any moving strand of elements in the ${\rm H}_{n}(q)$-module $\Lambda$ (\cite[Proposition~4]{DLP}). That is,
\[
X_{\widehat{h}}\, =\, X_{\widehat{bbm{(h)}}}\ \Leftrightarrow\ \underset{i}{\sum}\, \left(X_{\widehat{k_i}}\, =\, X_{\widehat{bbm{(k_i)}}}\right),
\]
\noindent where $h$ in the ${\rm H}_{n}(q)$-module $\Lambda^{aug}$ and $k_i$ in the ${\rm H}_{n}(q)$-module $\Lambda$.

\smallbreak

Finally, the `braiding tails' are eliminated from the elements in the ${\rm H}_{n}(q)$-module $\Lambda$. This reduces the computations to the set $\Lambda$, where the bbm's are executed on any moving strand (\cite[Theorem~8]{DLP}). Thus, in order to compute HOM$(L(p,1))$, it is sufficient to solve the infinite system of equations obtained by performing bbm's on every moving strand of elements in the set $\Lambda$. With a slight abuse of notation, we have that:
\[
HOM(L(p,1))\ =\ \frac{\Lambda}{<\lambda - bbm_i(\lambda)>},\ \lambda\in \Lambda, \ \forall\ i.
\]

\smallbreak

Moreover, in \cite{DL4} we consider the augmented set $\Lambda^{aug}$ and show that the system of equations obtained from elements in $\Lambda$ by performing bbm's on {\it any moving} strand, is equivalent to the system of equations obtained by performing {\it bbm}'s on the {\it first moving} strand of elements in $\Lambda^{aug}$. We show that by proving that the following diagram commutes:

\[
\begin{array}{ccccc}
\Lambda & \ni &  \tau_i & \overset{bbm_m}{\longrightarrow} & bbm_m(\tau_i)  \\
                 &     &     \updownarrow                &  &                              \updownarrow  \\
{\Lambda}^{aug} & \ni & \underset{j}{\sum}\tau_j & \overset{bbm_1}{\longrightarrow} &  \underset{j}{\sum} bbm_1(\tau_j)  \\
\end{array}
\]

\noindent Note that the fact that the $t_i$'s are not conjugate makes this procedure very non-trivial.

\bigbreak

The above are summarized in the following sequence of equations:

$$
\begin{array}{llllllll}
{\rm HOM}\left( L(p,1) \right) & = & \frac{{\rm HOM(ST)}}{<bbm's>} & = &
\frac{B_{1, n}}{<a - bbm_i(a)>, \ \forall\ i} & = & \frac{\Sigma_n^{\prime}}{<s^{\prime} - bbm_1(s^{\prime})>} & = \\
&&&&&&&\\
&&& = & \frac{\Sigma_n}{<s - bbm_1(s)>} & = & \frac{\Lambda^{aug}|{\rm H}_n}{<\lambda^{\prime} - bbm_i(\lambda^{\prime})>, \ \forall\ i} & = \\
&&&&&&&\\
&&& = & \frac{\Lambda|{\rm H}_n}{<\lambda^{\prime \prime} - bbm_i(\lambda^{\prime \prime})>, \ \forall\ i} & = & \frac{\Lambda}{<\lambda - bbm_i(\lambda)>, \ \forall\ i}& =\\
&&&&&&&\\
&&& = & \frac{\Lambda^{aug}}{<\mu - bbm_1(\mu)>}. &&&\\
\end{array}
$$

Namely:
\[
HOM(L(p,1))\ =\ \frac{\Lambda^{aug}}{<\mu - bbm_1(\mu)>},\ \mu\in \Lambda^{aug}.
\]

\begin{remark}\rm
It is worth mentioning that although $\Lambda^{aug} \supset \Lambda$, the advantage of considering elements in the augmented set $\Lambda^{aug}$ is that we restrict the performance of the braid band moves only on the first moving strand and, thus, we obtain less equations and more control on the infinite system. 
\end{remark}

Moreover, in \cite{D3} we consider the skein module obtained by solving the infinite system of equations $X_{\widehat{\mu}}\, =\, X_{\widehat{bbm_i{\mu}}}$, where we only consider elements in $\Lambda^{aug_+}$, a set related to $\mathcal{B}^+$, and we perform braid band moves on their first moving strands, namely $\frac{\Lambda^{aug_+}}{<bbm_1>}$. Similarly, we consider the module $\frac{\Lambda^{aug_-}}{<bbm_1>}$ obtained by solving the infinite system of equations by considering elements in $\Lambda^{aug_-}$, and we perform braid band moves on their first moving strands. These modules are related via two maps $f$ and $I$ defined as follows:

\begin{definition}\label{mapfI}\rm
\begin{itemize}
\item[(i)] We define the automorphism $f: \Lambda^{aug} \rightarrow \Lambda^{aug}$ such that:
\smallbreak
\[
\begin{array}{rcll}
f(\tau_1\cdot \tau_2) & = & f(\tau_1)\cdot f(\tau_2),& \forall\, \tau_1, \tau_2\in \Lambda^{aug}\\
&&&\\
t_i^k & \mapsto & t_{i}^{-k}, & \forall\ i\in \mathbb{N}^*,\ \forall\ k\in \mathbb{Z} \setminus \mathbb{N}\\\
\sigma_i & \mapsto & \sigma_i^{-1}, & \forall\ i\in \mathbb{N}^*
\end{array}
\]
\bigbreak
\item[(ii)] We define the map $I: R[z^{\pm 1}, s_k] \rightarrow R[z^{\pm 1}, s_k]$, $k\in \mathbb{Z}$ such that:
\smallbreak
\[
\begin{array}{rcll}
I(\tau_1 + \tau_2) & = & I(\tau_1) + I(\tau_2), & \forall\, \tau_1, \tau_2\\
I(\tau_1 \cdot \tau_2) & = & I(\tau_1) \cdot I(\tau_2), & \forall\, \tau_1, \tau_2\\
&&&\\
s_{-k} & \mapsto & s_{k}, & \forall\, k \in \mathbb{N}\\
s_{p-k} & \mapsto & s_{p+k}, & \forall\, k\ :\ 0 \leq k\leq p\\
z & \mapsto & \lambda \cdot z &\\
q^{\pm 1} & \mapsto & q^{\mp 1}&\\
\frac{\lambda^k}{z} & \mapsto & \frac{1}{\lambda^{k+1}z}, & \forall\, k\\
\end{array}
\]
\end{itemize}
\end{definition}

\begin{example}\rm
\begin{itemize}
\item[i.] $f(t^{-1}t_1^2\, \sigma_1^{-1})\, =\, tt_1^{-2}\, \sigma_1$.
\smallbreak
\item[ii.] For $p=3$ we have that $I(z\, q^{-2}\, s_{-2}\, s_2)\, =\, \lambda z\, q^{-2}\, s_2\, s_{4}$, since $I(s_2)\, =\, I(s_{3-1})\, :=\, s_{3+1}$.
\end{itemize}
\end{example}

Using the two maps $f$ and $I$ it is proved that the solution of the system $\frac{\Lambda^{aug_-}}{<bbm_1>}$ can be derived from the solution of the system $\frac{\Lambda^{aug_+}}{<bbm_1>}$. Indeed we have the following result:

\begin{theorem}\label{mthm2}
The equations obtained by imposing on the invariant $X$ relations coming from the performance of a $\pm $-$bbm_1$ on an element $\tau$ in $\Lambda^{aug_+}$ are equivalent to the image of the equations obtained by performing $\mp$-$bbm_1$ on its corresponding element $f(\tau)$ in $\Lambda^{aug_-}$ under $I$. That is:
$$I\left(X_{\widehat{f(\tau)}} =  X_{\widehat{bbm_{\mp 1}(f(\tau))}}\right)\ \Leftrightarrow \ X_{\widehat{\tau}} =  X_{\widehat{bbm_{\pm 1}(\tau)}}$$
\end{theorem}

Equivalently we have that the following diagram commutes:
\[
\begin{matrix}
\Lambda^{aug}_{(k)} & \ni & \tau & \overset{bbm_{\pm}}{\rightarrow} & bbm_{\pm 1}(\tau) & \Rightarrow & X_{\widehat{\tau}} =  X_{\widehat{bbm_{1}(\tau)}}, & X_{\widehat{\tau}} = X_{\widehat{bbm_{-1}(\tau)}}\\
                  &     & \updownarrow f &                     &                   &             &  \uparrow I &  \uparrow I\\
\Lambda^{aug}_{(-k)} & \ni & f(\tau) & \overset{bbm_{\mp}}{\rightarrow} & bbm_{\mp 1}(\tau) & \Rightarrow & X_{\widehat{f(\tau)}} =  X_{\widehat{bbm_{-1}(f(\tau))}}, & X_{\widehat{f(\tau)}} = X_{\widehat{bbm_{+1}(f(\tau))}}\\
\end{matrix}
\]

We now demonstrate Theorem~\ref{mthm2} with an example.

\begin{example}\rm
Consider the element $t\in \Lambda$ and perform a positive bbm to obtain $t^p\, t_1\, \sigma_1$. We have that
\[
\begin{array}{lclclc}
tr(t^p\, \underline{t_1}\, \sigma_1) & = & tr(t^p\, \sigma_1\, t\, \underline{\sigma_1^2}) & = & (q-1)\, tr(t^p\, \underline{t_1})\, +\, q\, tr(t^p\, \sigma_1\, \underline{t}) & =\\
&&&&&\\
& = & q(q-1)\, tr(t^p\, t_1^{\prime}) & + & (q-1)^2\, tr(t^{p}\, \sigma_1\, t)\, +\, q\, tr(t^p\, \sigma_1\, t) & =\\
&&&&&\\
& = & q(q-1)\, s_1s_p & + & (q^2-q+1)\, z\, s_{p+1} &
\end{array}
\]

\noindent Hence, we obtain the equation:
\[
\begin{array}{lclc}
X_{\widehat{t}}\, =\, X_{\widehat{t^pt_1\sigma_1}} & \Leftrightarrow & s_1\, =\, -\, \frac{1-\lambda q}{\sqrt{\lambda}(1-q)}\sqrt{\lambda}^3\cdot tr(t^pt_1\sigma_1) & \Leftrightarrow \\
&&&\\
&& s_1\, =\, \frac{\lambda}{z}\cdot \left[ q(q-1)\, s_1s_p \, + \, (q^2-q+1)\, z\, s_{p+1} \right]
\end{array}
\]

We now consider the homologous word of $t$, $t^{-1}:=f(t)$, and we perform a negative bbm, to obtain $bbm_{-1}(t^{-1})\, =\, t^p\, t_1^{-1}\, \sigma_1^{-1}$. We now evaluate $tr(t^p\, t_1^{-1}\, \sigma_1^{-1})$:
\[
\begin{array}{lclclc}
tr(t^p\, t_1^{-1}\, \underline{\sigma_1^{-1}}) & = & q^{-1}\, tr(t^p\, \underline{t_1^{-1}\, \sigma_1}) & + & (q^{-1}-1)\, tr(t^p\, \underline{t_1^{-1}}) & =\\
&&&&&\\
& = & q^{-1}\, tr(t^p\, \sigma_1^{-1}\, \underline{t^{-1}}) & + & (q^{-1}-1)\, tr(t^p\, \underline{\sigma_1^{-1}}\, t^{-1}\, \sigma_1^{-1}) & = \\
&&&&&\\
& = & q^{-1}\, tr(t^{p-1}\, \underline{\sigma_1^{-1}}) & + & (q^{-1}-1)\, q^{-1}\, tr(t^p\, {t_1^{\prime}}^{-1})\, +\, (q^{-1}-1)^2\, tr(t^{p-1}\, \underline{\sigma_1^{-1}}& = \\
&&&&&\\
& = & q^{-1}\, (q^{-1}-1)\, s_{-1}\, s_p & + & (\lambda\, z)\, (q^{-2}-q^{-1}+1)\, s_{p-1} &
\end{array}
\]

\noindent and we obtain the following equation:
\[
\begin{array}{lclc}
X_{\widehat{t^{-1}}}\, =\, X_{\widehat{t^pt_1^{-1}\sigma_1^{-1}}} & \Leftrightarrow & s_{-1}\, =\, -\, \frac{1-\lambda q}{\sqrt{\lambda}(1-q)}\sqrt{\lambda}^{-3}\cdot tr(t^pt_1^{-1}\sigma_1^{-1}) & \Leftrightarrow \\
&&&\\
&& s_{-1}\, =\, \frac{1}{\lambda^2 z}\cdot \left[ q^{-1}(q^{-1}-1)\, s_{-1}s_p \, + \, (q^{-2}-q^{-1}+1)\, \lambda\, z\, s_{p-1} \right]
\end{array}
\]

We now consider $I\left( X_{\widehat{t^{-1}}}\, =\, X_{\widehat{t^pt_1^{-1}\sigma_1^{-1}}} \right)$ and we have that
\[
s_1\, =\,I\left(\frac{1}{\lambda^2 z}\right)\cdot \left[ q(q-1)\, s_{1}s_p \, + \, (q^2-q+1)\, I\left(\lambda\, z\right)\, s_{p+1} \right],
\]

\noindent since $I(s_{-1})=s_1$, $I(s_{-1}s_p)=s_1s_p$, $I(s_{p-1})=s_{p+1}$ and $I(q^{-1})=q$. Finally, we have that
\[
I\left(\frac{1}{\lambda^2 z}\right)\, =\, I\left(\frac{\lambda^{-2}}{z}\right)\, :=\, \frac{\lambda}{z},\ {\rm and}\ I\left(\lambda\, z\right)\, =\, \lambda.
\] 

Hence, we obtain the equation $X_{\widehat{t}}\, =\, X_{\widehat{t^pt_1\sigma_1}}$.
\end{example}

Finally, in \cite{D3}, potential bases for $\frac{\Lambda^{aug_+}}{<bbm_1>}$ and for $\frac{\Lambda^{aug_-}}{<bbm_1>}$ are presented. In particular, it is shown that the set
\begin{equation}\label{setp}
\left\{{t^{\prime}}^{k_0}{t^{\prime}}^{k_1} \ldots {t_n^{\prime}}^{k_n},\ {\rm where}\ n,\, k_i\in \mathbb{N}\ :\ 0\leq k_i\leq p-1 \right\}
\end{equation}
\noindent is a generating set of the module $\frac{\Lambda^{aug_+}}{<bbm_1>}$, and since $\frac{\Lambda^{aug_+}}{<bbm_1>}\ =\ \frac{\Lambda^+}{<bbm_i>}$, we have that $\frac{\Lambda^+}{<bbm_i>}\ =\ \frac{{\Lambda^{\prime}}^+}{<bbm_i>}$. We conjecture that the set in Eq.~\ref{setp} forms a basis for $\mathcal{B}^+/<bbm_i>$. Similarly, the set
$$\left\{{t^{\prime}}^{k_0}{t^{\prime}}^{k_1} \ldots {t_n^{\prime}}^{k_n},\ {\rm where}\ n,\, k_i\in \mathbb{Z}\backslash \mathbb{N}\ :\ 0\geq k_i\geq -p+1 \right\}$$
\noindent is a generating set of the module $\frac{\Lambda^{aug_-}}{<bbm>}$. We conjecture that it also forms a basis for $\frac{\Lambda^{aug_-}}{<bbm>}$.

\begin{example}\rm
In this example we present some of the equations of the infinite system, a solution of which corresponds to the computation of HOM$(L(p,1))$:
\begin{itemize}
\item[i.] For $k=0$ we have $1 \in \Lambda^{aug_+}_{(0)}$ and $1\overset{bbm_{\pm }}{\rightarrow}t^p\sigma_1^{\pm 1}\, \Rightarrow\, 1:=s_{0}=s_{p}$.

\smallbreak

\item[ii.] For $k=1$ we have $t \in \Lambda^{aug_+}_{(1)}$ and $t\overset{bbm_{\pm }}{\rightarrow}t^pt_1\sigma_1^{\pm 1}\, \Rightarrow\, \begin{cases} s_{p+1}\, =s\, _{1}\\  s_ps_1\, =\, a_0s_1 \end{cases}$

\bigbreak

\item[iii.] For $k=2$ we have $ t^2, \, tt_1 \in \Lambda^{aug_+}_{(2)}$ and we obtain the following equations:
	\begin{itemize}
		\item[] $t^2 \overset{bbm}\rightarrow t^pt_1^{2}\sigma_1^{\pm 1}: \begin{cases} s_2 \ = \ a_1 s_{p+2} + a_2 s_{p+1}s_1 + a_3s_ps_2 \\
			s_2 \ = \ b_1s_{p+2} + b_2 s_{p+1}s_1 \end{cases}$
		
		\item[] $tt_1 \overset{bbm}\rightarrow t^pt_1t_{2}\sigma_1^{\pm 1}: \begin{cases} s_2 + s_1^2\ = \ c_1 s_{p+2} + c_2 s_{p+1}s_1 \\
		s_2 +s_1^2 \ = \ d_1 s_{p+2} + d_2 s_{p+1}s_1 + d_3s_ps_2 + d_4s_p s_1^2 \end{cases}$ 
			\end{itemize}
and thus:
\[\left\{\begin{array}{ccl}
			s_{p+2} & = & A_1 s_{2} + A_2 s_1^2 \\
			s_2s_p & = & B_1 s_{2} + B_2 s_1^2\\
            s_1s_{p+1} & = & C_1 s_{2} + C_2 s_1^2\\
            s_p s_1^2 & = & D_1 s_{2} + D_2 s_1^2\\
		\end{array}\right. \]
\bigbreak

\noindent where $a_i, A_i, b_i, B_i, c_i, C_i, d_i, D_i \in \mathbb{C},\ \forall\ i$.

\smallbreak

\item[iv.] For $k=-1$ we have $t^{-1} \in \Lambda^{aug_-}_{(-1)}$ and: $t^{-1}\overset{bbm_{\pm }}{\rightarrow}t^pt_1^{-1}\sigma_1^{\pm 1} \Leftrightarrow \begin{cases} s_{p-1}\, =\, s_{-1} \\ s_ps_{-1}\, =\, a^{\prime}_0s_1 \end{cases}$

\bigbreak

\item[v.] For the elements in $\Lambda^{aug_-}_{(-2)}$ we have:
	\begin{itemize}
		\item[] $t^{-2} \overset{bbm}\rightarrow t^pt_1^{-2}\sigma_1^{\pm 1}$
		\smallbreak
		\item[] $t^{-1}t_1^{-1} \overset{bbm}\rightarrow t^pt_1^{-1}t_{-2}\sigma_1^{\pm 1}$
and:
\[\left\{\begin{array}{ccl}
			s_{p-2} & = & A_1^{\prime} s_{-2} + A_2^{\prime} s_{-1}^2 \\
			s_{-2}s_p & = & B_1^{\prime} s_{-2} + B_2^{\prime} s_{-1}^2\\
            s_{-1}s_{p-1} & = & C_1^{\prime} s_{-2} + C_2^{\prime} s_{-1}^2\\
            s_p s_{-1}^2 & = & D_1^{\prime} s_{-2} + D_2^{\prime} s_{-1}^2\\
		\end{array}\right. \]
\bigbreak

\noindent where $A_i^{\prime}, B_i^{\prime}, C_i^{\prime}, D_i^{\prime} \in \mathbb{C},\ \forall\ i$.
\end{itemize}
\end{itemize}
\end{example}

\begin{remark}\rm
\begin{itemize}
\item[i.] Results on the infinite system so far suggest that the infinite system of equations admits unique solution, which on the level of skein modules means that the skein module is torsion free, and that the set
\[
\left\{{t^{\prime}}^{k_0}{t^{\prime}}^{k_1} \ldots {t_n^{\prime}}^{k_n},\ {\rm where}\ n,\, k_i\in \mathbb{Z}\ :\ -p/2 \leq k_i < p/2 \right\}
\]
\noindent forms a basis for HOM$(L(p, 1))$.
\smallbreak
\item[ii.] It is worth mentioning that this set was first presented in \cite{GM}, where, using diagrammatic methods, it is shown that it forms a basis for HOM$(L(p, 1))$. In \cite{DL5} we work toward proving this result using braids and techniques described within this paper.
\smallbreak
\item[iii.] We have reasons to believe that the braid technique can be successfully applied in order to compute skein modules of other more complicated c.c.o. 3-manifolds (such as the lens spaces $L(p, q), q>1$), where diagrammatic methods fail to do so (\cite{DGLM}).
\end{itemize}
\end{remark}

\subsection{The Kauffman bracket skein module of the lens spaces $L(p,1)$ for $p\neq 0$}\label{kbsmlp1}

In this subsection we present a new basis for KBSM($L(p,1)$), $p\neq 0$ first appeared in \cite{D0} (see also \cite{D1}). This basis is different from that presented in \cite{HP} and it so more natural on the level of braids, since elements in this basis have no crossings.

\smallbreak

In \cite[Theorem 3]{D0} the following result is proved:

\begin{theorem}\label{mthm11}
The set $\mathcal{B}_p\ := \{t^n,\, n\in I \}$, where $I=\{0, 1, \ldots, \lfloor p/2 \rfloor\}$, is a basis for the Kauffman bracket skein module of the lens spaces $L(p,1)$.
\end{theorem}

In order to prove Theorem~\ref{mthm11} we first translate the relation ${\rm KBSM}(L(p,1))\, =\, \frac{{\rm KBSM(ST)}}{<a-bm(a)>}$, for all $a$ in a basis of KBSM(ST), in terms of braids. Namely, we impose on the universal invariant $V$ of the Kauffman bracket type for knots and links in ST, relations coming from the performance of bbm's on elements in the $B_{\rm ST}$ basis of KBMS(ST). That is, in order to compute KBSM($L(p,1)$) we solve the infinite system of equations: 

\[
V_{\widehat{t^n}}\, =\, V_{\widehat{t^pt_1^n\sigma_1^{\pm 1}}}
\] 

Note that the unknowns of the infinite system are the $s_i$'s that come from the fourth rule of the Markov trace.

\smallbreak

For computing $tr\left(t^pt_1^n\sigma_1^{\pm 1}\right)$, we first express $t^pt_1^n\sigma_1^{\pm 1}$ to sums of elements of the form $t^m{t_1^{\prime}}^k$. In particular, the following formulas hold in $TL_{1, n}$ for $n, p \in \mathbb{N}\backslash \{0\}$:

\[
t^p\, t_1^n\, \sigma_1\ \widehat{\cong}\ \underset{i=0}{\overset{n}{\sum}}\, a_i\, t^{p+i}\, {t_1^{\prime}}^{n-i}\ \quad\ {\rm and}\ \quad\ t^p\, t_1^n\, \sigma_1^{-1}\ \widehat{\cong}\ \underset{i=1}{\overset{n}{\sum}}\, b_i\, t^{p+i}\, {t_1^{\prime}}^{n-i}
\]

\noindent where $a_i, b_i$ coefficients for all $i$ and the symbol $\widehat{\cong}$ denotes that conjugation and stabilization moves are performed. Finally, we translate the resulting elements in terms of elements in the basis $B_{\rm ST}$ and we have that:

\begin{equation}\label{l2}
t^p\, {t_1^{\prime}}^n \ \widehat{\underset{{\rm skein}}{\cong}}\ \underset{i=0}{\overset{\lfloor n/2\rfloor}{\sum}}\, a_i \, t^{p+n-2i},
\end{equation}
\noindent where $a_i$ are coefficients for all $i$.

\smallbreak

Hence, we have that $tr(t^p\, t_1^n\, \sigma_1)\ =\ \underset{i=0}{\overset{\lfloor n/2\rfloor}{\sum}}\, a_i \, s_{p+n-2i}$ for some coefficients $a_i$, since:
\[
tr(t^p\, t_1^n\, \sigma_1)\ =\ tr\left(\underset{i=0}{\overset{n}{\sum}}\, a_i\, t^{p+i}\, {t_1^{\prime}}^{n-i}\right)\ =\ tr\left(\underset{i=0}{\overset{n}{\sum}}\, \underset{j=0}{\overset{\lfloor \frac{n-i}{2}\rfloor}{\sum}}\, \left(a_i^{\prime}\, t^{p+n-2j}\right)\right).
\]

\noindent The following diagram summarizes the steps we followed in order to obtain the equations of the infinite system:

\[
\begin{array}{ccccccc}
t^n & \overset{bbm}{\longrightarrow} & t^pt_1^n\, \sigma_1 & \hat{\cong} & \underset{i=0}{\overset{n}{\sum}}\, a_i\, t^{p+i}\, {t_1^{\prime}}^{p-i} & \hat{\underset{{\rm skein}}{\cong}} & \underset{i=0}{\overset{n}{\sum}}\, \underset{j=0}{\overset{\lfloor \frac{n-i}{2} \rfloor}{\sum}}\, \left(s_{p+n-2j}\right) \\
|   &                 &                     &                 &                                                                          &                                        &       |\\
tr   &                 &                     &                 &                                                                          &                                        &       tr\\
\downarrow   &                 &                     &                 &                                                                          &                                        &       \downarrow\\
s_n &                 &                    &                  &                                                                           &                                                                                                                  &        \underset{i=0}{\overset{n}{\sum}}\, s_{p-n+2i}          \\      
|   &                 &                     &                 &                                                                          &                                        &       |\\
V   &                 &                     &                 &                                                                          &                                        &       V\\
\downarrow   &                 &                     &                 &                                                                          &                                        &       \downarrow\\   
s_n &                 &                    &                  &            =                                                               &                                                                                                                  &        \left(-\, \frac{1+u^2}{u} \right)\, u^{2e}\cdot \underset{i=0}{\overset{n}{\sum}}\, s_{p-n+2i}          \\
\end{array}
\]

We now treat elements of the form $s_{p+k}$ for $k\in \mathbb{N}$ and we show that these elements can be written as sums of elements in $s_i$'s with $i\in \{0, 1, \ldots, p\}$. In particular, in \cite[Proposition 3]{D0} the following formulas are proved:

\[
s_{p+n}\ =\ \begin{cases} \underset{i=0}{\overset{n/2}{\sum}}\, (c_i\, s_{2i}\, +\, d_i\, s_{p-2i}) &,\, {\rm for}\ n\ {\rm even}\\
\underset{i=0}{\overset{(n-1)/2}{\sum}}\, (c^{\prime}_i\, s_{2i-1}\, +\, d^{\prime}_i\, s_{p-2i-1}) &,\,{\rm for}\ n\ {\rm odd} \end{cases},
\]
\noindent where $c_i, d_i, c^{\prime}_i, d^{\prime}_i$ are coefficients. Namely, for $p, k\in \mathbb{N}$ we have that:
\[
s_{p+k}\ =\ \underset{i<p}{\sum}\, q_i\, s_i, \quad {\rm where}\ i\in \mathbb{Z}\ {\rm and}\ q_i\ {\rm coefficients}.
\]

\smallbreak

Note now that the indices on some monomials of the resulting sum may be negative. As shown in \cite{D0}, these monomials may be expressed in terms of $s_i$'s where $i\in \mathbb{N}$. In particular we have that:

\[
s_{-n}\ = \begin{cases} \underset{i=0}{\overset{n/2}{\sum}}\, s_{2i} &, {\rm for}\ n\ {even}\\
\underset{i=0}{\overset{(n-1)/2}{\sum}}\, s_{2i+1} &, {\rm for}\ n\ {odd}
\end{cases}
\]

Hence, we conclude that the unknowns of the infinite system of equations with indices greater than or equal to $p$ can be written as sums of unknowns with non-negative indices and also that these indices are $\mod p$. In other words, we conclude that the set $\mathcal{B}_p$ spans the Kauffman bracket skein module of the lens spaces $L(p,1)$. Finally \S~4 \cite{D0}, the set $\mathcal{B}_p$ is shown to be linearly independent and hence it forms a basis for KBSM($L(p,1)$).

\begin{example}\rm
In this example we present equations of the infinite system and we demonstrate the above mentioned results. We have that:  
\[
\begin{array}{lclcl}
 t^0\ \overset{bbm}{\rightarrow}\ t^p\, \sigma_1^{\pm 1} & \overset{V}{\Rightarrow} & 1:=s_0\, =\, s_p & \Rightarrow & s_p\ =\ s_0\\
&&&&\\
t^{-1}\ \overset{bbm}{\rightarrow}\ t^p\, t_1^{-1}\sigma_1^{\pm 1} & \overset{V}{\Rightarrow} & s_{-1}\, =\, s_{p-1} & {\underset{s_{-1}=s_1}{\Rightarrow}}& s_{p-1}\ =\ s_1\\
&&&&\\
t^{-2}\ \overset{bbm}{\rightarrow}\ t^p\, t_1^{-2} \sigma_1^{\pm 1} & \overset{V}{\Rightarrow} & s_{-2}\, =\, s_{p-2}+s_p+s_{p+2} & {\underset{\left\{\begin{matrix} s_{-2}& = & s_0+s_2\\ s_{p+2} & = & \underset{i=0}{\overset{p}{\sum}}\, s_i\end{matrix}\right\}}{\Rightarrow}}& s_{p-2}\ =\ \underset{i=0}{\overset{p-3}{\sum}}\, s_i\\
&&&&\\
 &\vdots& &\vdots& \\
\end{array}
\]

\[
\begin{array}{lcl}
t^{-\lfloor p/2 \rfloor +1}\ \overset{bbm}{\rightarrow}\ t^p\, t_1^{-\lfloor p/2 \rfloor +1} \sigma_1^{\pm 1} & \overset{V}{\Rightarrow} & s_{-\lfloor p/2 \rfloor +1}\, =\, s_{p-\lfloor p/2 \rfloor -1}+\ldots+s_{p+\lfloor p/2 \rfloor+1}\\
&&\\
& \Rightarrow & s_{p-\lfloor p/2 \rfloor +1}\ =\ \underset{i=0}{\overset{p-\lfloor p/2 \rfloor}{\sum}}\, s_i
\end{array}
\]

\end{example}

\section{The Kauffman bracket skein module of $S^1\times S^2$}\label{kbsms1s2}

\begin{figure}
\begin{center}
\includegraphics[width=2.9in]{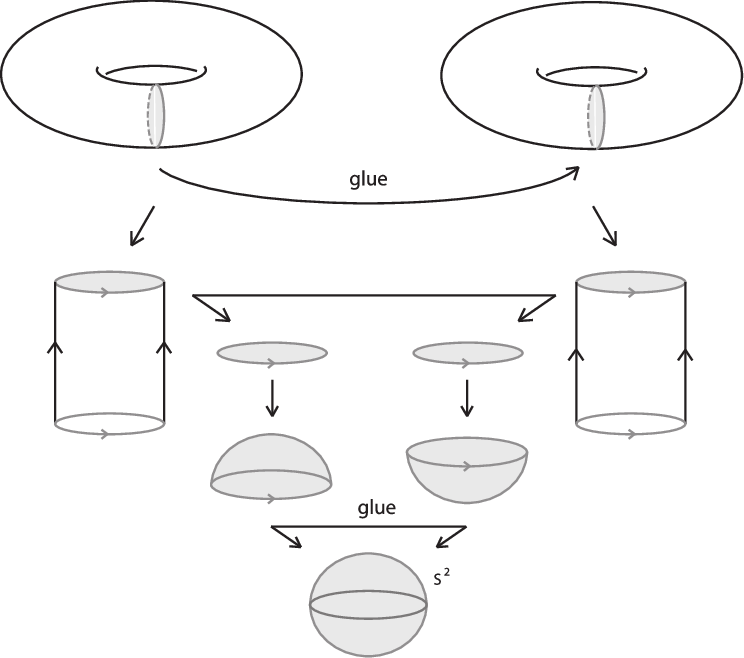}
\end{center}
\caption{Gluing two solid tori to obtain $S^1 \times S^2$.}
\label{stst}
\end{figure}

$S^1\times S^2$ may be obtained by `gluing' two solid tori via some homeomorphism $h$ on their boundaries, such that $h$ takes a meridian curve $m_1$ on the first solid torus ${\rm ST}_1$, to the corresponding meridian curve $m_2$ of the second solid torus ${\rm ST}_2$. We may choose $h$ to be the identity $i: \partial {\rm ST}_1\, \rightarrow\, \partial {\rm ST}_2$. Since now $ST\, =\, S^1\times D^2$, where $D^2$ denotes a disc, and since the gluing of two discs results in a 2-sphere $S^2$, the resulting $3$-manifold can be considered as a family of $2$-spheres parametrized by a circle (for an illustration see Figure~\ref{stst}). 

In this subsection we compute the Kauffman bracket skein module of $S^1\times S^2$, ${\rm KBSM}\left(S^1\times S^2\right)$, via braids. We follow the same procedure as in \S~\ref{kbsmlp1} and we extend the universal Kauffman bracket type invariant $V$ for knots and links in the Solid Torus to an invariant for knots and links in $S^1\times S^2$. In particular, we show that ${\rm KBSM}\left(S^1\times S^2\right)$ is not torsion free and that its free part is generated by the unknot (or the empty knot).

\subsection{Topological set up}\label{setup2}

\begin{figure}[H]
\begin{center}
\includegraphics[width=3.9in]{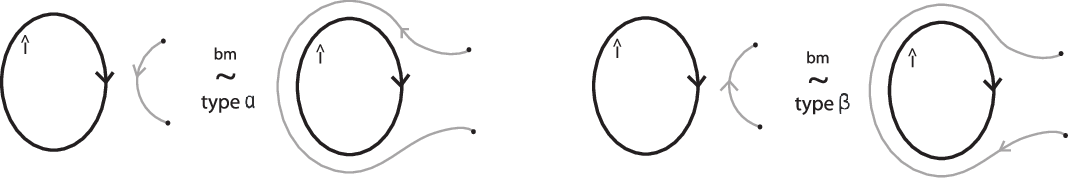}
\end{center}
\caption{The two types of band moves for links in $S^1 \times S^2$.}
\label{bmos1s2}
\end{figure}

We consider $S^1\times S^2$ as being obtained from $S^3$ by integral surgery along the unknot with coefficient zero, that is, $S^1\times S^2$ is a special type of lens spaces $L(p, q)$ for $p=0$ and $q=1$ ($L(0,1)\, \cong\, S^1\times S^2$). The band moves, that reflect isotopy in $S^1 \times S^2$ are illustrated in Figure~\ref{bmos1s2}.

\begin{figure}[H]
\begin{center}
\includegraphics[width=2.8in]{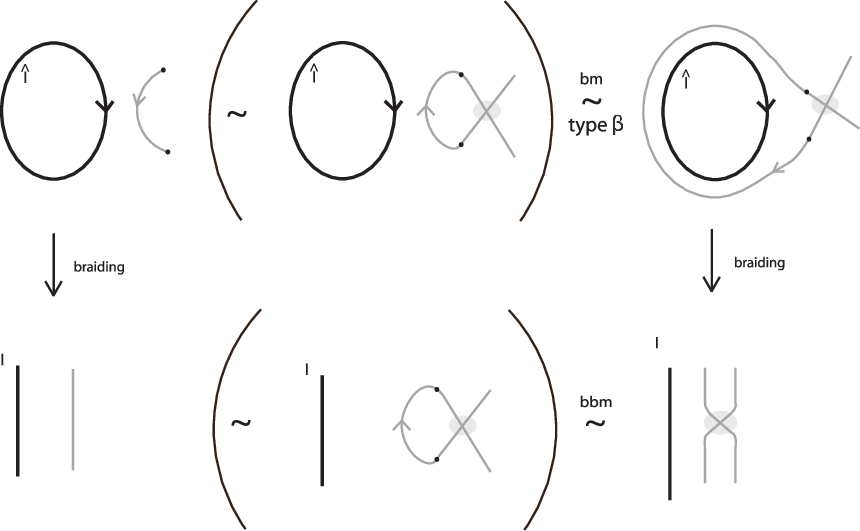}
\end{center}
\caption{ The two types of braid band moves for braids in $S^1 \times S^2$.}
\label{bbmovs1s2}
\end{figure}

Recall that a braid band move is defined as an equivalence move on mixed braids that corresponds to a band move on th closure of the mixed braid. In the case of $S^1 \times S^2$, these moves have the following algebraic expression:
\[
\alpha \sim  \alpha_+ \sigma_1^{\pm 1},
\]
\noindent where $\alpha\in B_{1,n}$ and $\alpha_+\in B_{1, n+1}$ is the word $\alpha$ with all indices shifted by +1. For an illustration see Figures~\ref{bbmovs1s2} and \ref{tbbm1}.

\subsection{The Kauffman bracket skein module of $S^1\times S^2$}

\begin{figure}[H]
    \includegraphics[width=1.7in]{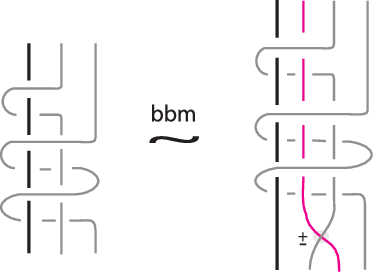}
    \caption{Braid band moves performed on $tt_1^2 \in B_{ST}$.}
    \label{tbbm1}
\end{figure}

In this subsection we solve the infinite system of equations obtained by performing braid band moves on elements in the basis $B_{\rm ST}$ of KBSM(ST) and by imposing to the generic invariant $V$ for knots and links in ST relations of the form $V_{\widehat{t^n}}\, =\, V_{\widehat{bbm(t^n)}}$, for all $n\in \mathbb{N}$, where $bbm(t^n)\, =\, t_1^n\, \sigma_1^{\pm 1}$ (see Figure~\ref{tbbm1}). Recall that the unknowns in the system are the $s_i$'s, coming from the fourth rule of the trace function in Theorem~\ref{tr}, that is, $tr(t^n)\, =\, s_n$, for all $n\in \mathbb{Z}$.

The method used for computing KBSM$(S^1 \times S^2)$ differs significantly from the one used in computing KBSM$(L(p,1))$. More precisely, we start by considering the effect of negative bbm's on elements in $B_{\rm ST}$ and we obtain a spanning set for $B^{S^1\times S^2}_{-}\, :=\, \frac{B_{\rm ST}}{<\, a-bbm_{-}(a)\, >}$. Then we study the effect of positive bbm's on elements in $B^{S^1\times S^2}_{-}$. With a slight abuse of notation, we denote the above as follows:

\[
{\rm KBSM}\left(S^1 \times S^2\right)\, =\, \frac{{\rm B}_{ST}}{<\, a-bbm_{\pm}(a)\, >}\, =\, \frac{{\rm B}_{ST}/<a-bbm_-(a)>}{<a-bbm_+(a)>}\, =\, \frac{B^{S^1\times S^2}_{-}}{<a-bbm_+(a)>}.
\]

Hence, we consider elements in $B_{ST}$ and perform negative braid band moves, obtaining an infinite system of equations, the solution of which corresponds to a basis for $B^{-}_{S^1\times S^2}$. 

\begin{example}\rm
In this example we show that $t^0$ and $t$ are free in $B^{-}_{S^1\times S^2}$. 
\[
\begin{array}{rcrcl}
1 \ \overset{bbm_-}{\rightarrow}\ \sigma_1^{-1} & \Rightarrow & V_{\hat{1}}\, =\, V_{\hat{\sigma_1^{-1}}} & \Rightarrow & 1 \ = \ \left(-\, \frac{1+u^2}{u} \right)\, u^{-2}\, tr(\sigma_1^{-1})\, \Rightarrow\\
&&&&\\
 & & 1 & = & -\, \frac{1+u^2}{u^3}\, \left(z-(u-u^{-1}) \right)\ \overset{z=-\, \frac{1}{u(1+u^2)}}{\Rightarrow}\\
&&&&\\
 & & 1 & = & 1 \\
&&&&\\
 t \ \overset{bbm_-}{\rightarrow}\ t_1\, \sigma_1^{-1} & \Rightarrow & V_{\hat{t}}\, =\, V_{\widehat{t_1\, \sigma_1^{-1}}} & \Rightarrow & s_1 \ = \ \left(-\, \frac{1+u^2}{u} \right)\, u^{2}\, tr(t_1\, \sigma_1^{-1})\, \Rightarrow\\
&&&&\\
 & & s_1 & = & -\, -\, \frac{1+u^2}{u}\, u^2\, z\, s_1\ \Rightarrow\\
&&&&\\
 & & s_1 & = & s_1\\
\end{array}
\]
\end{example}

In \cite{D10}, the following result is presented:
\[
\begin{array}{lcc}
V_{\hat{t^n}}\, =\, V_{\widehat{t_1^n\sigma_1^{-1}}} & \Rightarrow & s_n\, =\, \begin{cases} \underset{i=0}{\overset{(n-2)/2}{\sum}}\, s_{2i} &,\ {\rm for}\ n\ {\rm even}\\
\underset{i=0}{\overset{(n-3)/2}{\sum}}\, s_{2i+1} &,\ {\rm for}\ n\ {\rm odd}
\end{cases}
\end{array}
\]

Equivalently, the elements $t^n \in B_{\rm ST}, n\geq 2$ can be written as sums of $t^0$ and $t$. This leads to the following theorem:

\begin{theorem}\label{thm12}
$B^{-}_{S^1\times S^2}$ is generated by the unknot $t^0$ and $t$. 
\end{theorem}

From Theorem~\ref{thm12} and from the fact that the elements $t^0$ and $t$ are linearly independent, we have that the set $\{t^0,\, t\}$ forms a basis for $B^{-}_{S^1\times S^2}$. We then consider the effect of positive braid band moves on the elements $t^0$ and $t$. We have that:

\[
\begin{array}{lrcrcl}
{\rm For}\ t^0: & 1 \ \overset{bbm_+}{\rightarrow}\ \sigma_1 & \Rightarrow & V_{\hat{1}}\, =\, V_{\hat{\sigma_1}} & \Rightarrow & 1 \ = \ \left(-\, \frac{1+u^2}{u} \right)\, u^{-2}\, tr(\sigma_1)\, \Rightarrow\\
&&&&&\\
& & & 1 & = & -\, (1+u^2)\, u\, z \ \overset{z=-\, \frac{1}{u(1+u^2)}}{\Rightarrow}\\
&&&&&\\
& & & 1 & = & 1 \qquad {\rm (free\ part)}
\end{array}
\] 

Hence, the unknot is free in KBSM($S^1\times S^2$. We now deal with the generator $t$ and we have the following:
\[
\begin{array}{rcrcll}
t \ \overset{bbm_+}{\rightarrow}\ t_1\, \sigma_1 & \Rightarrow & V_{\hat{t}}\, =\, V_{\widehat{t_1\, \sigma_1}} & \Rightarrow & s_1 \ = \ \left(-\, \frac{1+u^2}{u} \right)\, u^{6}\, tr(t_1\, \sigma_1) & \overset{tr(t_1\, \sigma_1)=(u^2-1+u^{-2})\, z\, s_1\ +\ (u-u^{-1})\, s_1}{\Rightarrow}\\
&&&&&\\
&&&\Rightarrow & (1-u^6)\, (1-u^2)\, s_1 \ =\  0. & \\
\end{array}
\]

The discussion above leads to the following theorem:

\begin{theorem}\label{mthm12}
The free part of KBSM$\left(S^1 \times S^2\right)$ is generated by the unknot and the torsion part is generated by $t$.
\end{theorem}

Equivalently, we have that
\[
{\rm KBSM}\left(S^1 \times S^2\right)\ =\ \mathbb{Z}[A^{\pm 1}]\, \oplus\, {\rm Torsion}
\]

\noindent and although we don't obtain a closed formula for the torsion part of KBSM($S^1 \times S^2$), Theorem~\ref{mthm12} implies the following:
\[
\pi(t^{2n+1})\, =\, 0,\, \forall\, n\in \mathbb{N} \quad {\rm and} \quad {\rm KBSM}\left(S^1 \times S^2\right)/Tor\, =\, \mathbb{Z}[A^{\pm 1}],
\]
\noindent where $\pi$ denotes the natural projection $\pi: {\rm KBSM}\left(S^1 \times S^2\right)\, \rightarrow\, {\rm KBSM}\left(S^1 \times S^2\right)/Tor$.

\section{Conclusions}

In this paper we present recent results on the computation of the HOMFLYPT and the Kauffman bracket skein modules of $S^3$, the Solid Torus, and the lens spaces. We first present new bases for HOM(ST) and KBSM(ST) in braid form and we then relate HOM$(L(p,1))$ to HOM(ST), and KBSM$(L(p,1))$ to KBSM(ST), by means of equations resulting from the performance of {\it bbm}'s. For both skein modules, we arrive at infinite systems of equations, the solution of which coincides with the computation of the corresponding skein module of $L(p,1)$.

\smallbreak

For the case of KBSM, we solve this system and we show that the set

\[
\mathcal{B}_p\ := \{t^n,\, n\in I \},\ {\rm where}\ I=\{0, 1, \ldots, \lfloor p/2 \rfloor\}
\]

\noindent forms a basis for the Kauffman bracket skein module of the lens spaces $L(p,1)$ (Theorem~\ref{mthm11}). For KBSM$\left(S^1 \times S^2\right)$, we show that it is generated by the unknot and the torsion part is generated by $t$ (Theorem~\ref{mthm12}).

\smallbreak

For the case of the HOMFLYPT skein module of the lens spaces $L(p,1)$, we show that in order to compute HOM$(L(p,1))$ it suffices to solve the infinite system of equations:

$$X_{\widehat{\tau}}\ =\ X_{\widehat{bbm_i(\tau)}},$$

\noindent where $bbm_i(\tau)$ is the result of the performance of {\it bbm} on the $i^{th}$-moving strand of $\tau \in \Lambda$, for all $\tau \in \Lambda$ and for all $i$.

\noindent and we conclude that the set
\[
\left\{{t^{\prime}}^{k_0}{t^{\prime}}^{k_1} \ldots {t_n^{\prime}}^{k_n},\ {\rm where}\ n,\, k_i\in \mathbb{Z}\ :\ -p/2 \leq k_i < p/2 \right\}
\]

\noindent forms a potential basis for HOM$(L(p,1))$.

\smallbreak

As a final note, it is worth mentioning that a different diagrammatic method based on braids has been successful in the computation of Kauffman bracket skein modules of 3-manifolds. More precisely, this method has been successfully applied for the case of the lens spaces $L(p,q)$ in \cite{D0}, for $S^1\times S^2$ in \cite{D10}, for the handlebody of genus 2 in \cite{D2} and for the complement of $(2,2p+1)$-torus knots in \cite{D4}.

\end{document}